\documentclass[reqno]{amsart}
\usepackage{amssymb}
\usepackage{amsmath}
\usepackage{changes}
\usepackage[mathscr]{euscript}
\usepackage[small]{caption}
\usepackage{mathtools}
\usepackage{graphicx}
\graphicspath{ {images/} }
\usepackage{xcolor}
\usepackage{comment}
\usepackage{changes}

\makeatletter
\@addtoreset{equation}{section}
\makeatother

\newtheorem{theorem}{Theorem}[section]

\newtheorem{remark}[theorem]{Remark}
\newtheorem{example}[theorem]{Example}

\newcounter{as}[section]

\newcommand{\mc}[1]{{\mathcal #1}}
\newcommand{\mf}[1]{{\mathfrak #1}}

\newcommand{\bb}[1]{{\mathbb #1}}
\newcommand{\bs}[1]{{\boldsymbol #1}}
\newcommand{\ms}[1]{{\mathscr #1}}

\newcommand{\<}{\langle}
\renewcommand{\>}{\rangle}

\definecolor{bblue}{rgb}{.2,0.2,.8}

\title[Systems in
mild contact with boundary reservoirs]
{Thermodynamics of nonequilibrium driven diffusive systems in
mild contact with boundary reservoirs.
}

\author[A. Bouley]{Ang\`ele Bouley} 
\address{Ang\`ele Bouley
  \hfill\break\indent CNRS UMR 6085, Universit\'e de
  Rouen, \hfill\break\indent Avenue de l'Universit\'e, BP.12,
  Technop\^ole du Madril\-let, \hfill\break\indent
F76801 Saint-\'Etienne-du-Rouvray, France.} 
\email{angele.bouley@univ-rouen.fr}

\author[C. Landim]{Claudio Landim} 
\address{Claudio Landim
  \hfill\break\indent IMPA \hfill\break\indent Estrada Dona Castorina
  110, \hfill\break\indent
J. Botanico, 22460 Rio de Janeiro, Brazil\hfill\break\indent
  {\normalfont and} \hfill\break\indent CNRS UMR 6085, Universit\'e de
  Rouen, \hfill\break\indent Avenue de l'Universit\'e, BP.12,
  Technop\^ole du Madril\-let, \hfill\break\indent
F76801 Saint-\'Etienne-du-Rouvray, France.} 
\email{landim@impa.br}

\begin{document}

\begin{abstract}
We consider macroscopic systems in mild contact with boundary
reservoirs and under the action of external fields. We present an
explicit formula for the Hamiltonian of such systems, from which we
deduce the equation of motions, the action functional, the
hydrodynamic equation for the adjoint dynamics, and a formula for the
quasi-potential.

We examine the case in which the external forcing depends on time and
drives the system from one nonequilibrium state to another. We extend
the results presented in \cite{bgjl1} on thermodynamic transformations
for systems in strong contact with boundary reservoirs to the present
situation.

In particular, we propose a natural definition of renormalized work,
and show that it satisfies a Clausius inequality, and that quasi-static
transformations minimize the renormalized work. In addition, we
connect the renormalized work to the quasi-potential describing the
fluctuations in the stationary nonequilibrium ensemble.
\end{abstract}


\maketitle

\section{Introduction}

After the recent articles \cite{DHS, BEL21}, where a formula for the
quasi-potential has been derived for one-dimensional exclusion
processes in mild contact with reservoirs, the purpose of this article
is to extend to driven diffusive systems in mild contact with boundary
reservoirs the nonequilibrium thermodynamical theory developed in
\cite{bgjl1, bgjl2, bdgjl14} for systems in strong interaction with
reservoirs.

The macroscopic evolution of systems in mild contact with reservoirs
differs substantially from the one observed when the system strongly
interacts with the reservoirs. With strong interactions, at the level
of large deviations, a density fluctuation at the boundary is too
costly and not observed. In contrast, for mild boundary interactions,
any smooth trajectory has a finite cost. In consequence, the
associated Hamiltonian carries a term which takes into account the
boundary fluctuations.  We investigate in this article the
consequences to the thermodynamical theory of the additional boundary
Hamiltonian term.

In Section \ref{sec1}, we introduce a class of stochastic lattice
gases, which includes exclusion, zero-range and KMP models.  We define
in this framework a boundary Hamiltonian, derived rigorously in
\cite{FGLN2021, BEL21} for one-dimensional, symmetric exclusion
processes in mild contact with reservoirs, and record some of its
properties. For the readers convenience, these dynamics are reviewed
in the appendices, where explicit formulas for the Hamiltonians and
the quasi-potentials are presented.  We also provide a microscopic
dynamics not covered by this theory, as the stationary states of the
boundary dynamics are different from the bulk ones.

In Section \ref{sec0}, supported by the models presented in the
previous section, we introduce the main object of this article, a
Hamiltonian composed of a bulk part and the boundary part already
put forward.  The bulk part coincides with the Hamiltonian of systems in
strong interaction with the reservoirs, and is expressed in terms of
two thermodynamics quantities, the diffusivity and the mobility. The
boundary Hamiltonian, instead, the main novelty of this article, is
expressed in terms of a measure which depends on the chemical
potential of the reservoirs and on the density profile at the
boundary. We present in Section \ref{sec0} several properties of this
Hamiltonian and deduce from its form the equation of motions, the
action functional, the quasi-potential, the adjoint thermodynamical
quantities and formulas for the currents.

In Section \ref{sec6}, we derive a differential equation for the
quasi-potential. This equation has been obtained by Derrida,
Hirschberg and Sadhu \cite{DHS} for one-dimensional symmetric
exclusion processes, expressing the stationary state of the system as
a product of matrices. For zero-range processes it can be obtained by
direct computations because the stationary state is a product
measure. For KMP models, the equation is new and has not yet been
derived rigorously. In this model, the boundary conditions, displayed
in equations \eqref{3-13d} and \eqref{3-13db}, do not coincide with
the boundary conditions for stationary density profile (see equation
\eqref{3-09}), in contrast to the case of strong boundary
interactions and of exclusion and zero-range dynamics in mild
interaction. See Remarks \ref{rm-x} and \ref{rm-x2}.

Section \ref{sec3} provides a dynamic derivation of the second law of
thermodynamics as expressed by a Clausius inequality for the energy
exchanged between the system and the external reservoirs and fields.
The results and the reasoning presented in this section and the next
one follow closely \cite{bgjl1, bgjl2}.

In Section \ref{sec4}, we examine transformations along equilibrium
states. According to the standard thermodynamic theory, a
transformation is reversible if the energy exchanged between the
system and the environment is minimal. A thermodynamic principle
asserts that reversible transformations are accomplished by a sequence
of equilibrium states and are well approximated by quasi-static
transformations in which the variations of the environment are very
slow. By an explicit construction of quasi-static transformations, we
show that this principle can be derived for driven diffusive systems
in mild interaction with boundary reservoirs.

Fix two equilibrium states and a transformation which drives the
system from the first to the second one.  The \emph{excess work} of
this transformation is defined as the total work minus the minimal
work needed to bring the system from the first to the second
equilibrium state.  We show in Section \ref{sec4} that the
quasi-potential coincides with the excess work of the relaxation path
from the first equilibrium state to the second.

In Section \ref{sec5}, these results are extended to transformations
along nonequilibrium states. However, nonequilibrium states are
characterized by the presence of a non vanishing current in the
stationary density profile. Therefore, to maintain such states one
needs to dissipate a positive amount of energy per unit of time.  If
we consider a transformation between nonequilibrium stationary states,
the energy dissipated along such transformation will necessarily
include the contribution needed to maintain such states.

To take into account this amount of energy, following \cite{op}, we
introduce the \emph{renormalized work}, and extend the results of the
two previous sections to transformations along nonequilibrium
states. In contrast with systems in strong interaction with the
reservoirs, the symmetric and the anti-symmetric currents are not
orthogonal (cf. the discussion at the end of Section \ref{sec0}). In
consequence, the definition of the renormalized work proposed here is
new and involves a boundary functional which takes into account the
interaction of the system with the reservoirs.

The Hamiltonian appearing in this article is the action functional of
the dynamical large deviations principle (DLDP) for the empirical
measure.  A dynamical large deviations principle for systems in mild
contact with reservoirs has only been derived rigorously for
one-dimensional symmetric exclusion processes \cite{FGLN2021, BEL21}.
We believe that a DLDP for systems in mild contact with reservoirs can
also be derived for gradient exclusion processes in any dimension. For
zero-range and KMP processes, however, a rigorous proof is still out
of reach due to a lack of exponential moments \cite{BGL05}.

A formula for the Hamiltonian, as explained in Section \ref{sec1}
(cf. equation \eqref{3-02} and \eqref{x11}), can be easily obtained 
from the generator of the boundary dynamics and the stationary state
of the bulk dynamics.

A large deviations principle for the empirical measure under
nonequilibrium states for systems in mild contact with reservoirs is
more demanding and has not yet been proved. A formula for the rate
functional (the quasi-potential) for one-dimensional symmetric
exclusion processes is presented in \cite{DHS, BEL21}. For
zero-range processes, as the nonequilibrium states are product
measures, it is easy to derive it. For all other models, it is still
an open problem.

In conclusion, in this article we extend the thermodynamic theory
resulting from the Macroscopic Fluctuation Theory (MFT) to systems in
mild contact with reservoirs. We introduce a Hamiltonian with an
additional term coming from the mild interactions of the system with
reservoirs and propose a definition of renormalized work. We prove the
validity of a Clausius inequality for transformations along
equilibrium and nonequilibrium states in this framework and we show
that the excess of work along the relaxation path is given by the
quasi-potential, for equilibrium and nonequilibrium states.

\section{Microscopic dynamics}
\label{sec1}

In this section, we introduce the boundary Hamiltonian from an
underlying microscopic dynamics. The evolution induced by this
Hamiltonian together with the bulk Hamiltonian arising from locally
conservative dynamics will be examined in the next sections.

The general framework presented in this section includes the main
microscopic stochastic dynamics, such as the exclusion, zero-range and
KMP processes, on which the Macroscopic Fluctuation Theory (MFT) has
been build. For the reader's convenience, we reviewed
in the appendices the properties of these systems used below.

Let $\Omega$ be the bounded domain of $\bb R^d$ occupied by the
system. Fix $N\ge 1$, and denote by
$\color{blue} \Omega_N = \Omega \cap (\bb Z/N)^d$ its
discretization. Here, $\bb Z/N = \{k/N: k\in\bb Z\}$.  Elements of
$\Omega$ are represented by $x$, $y$.  The boundary of $\Omega_N$,
denoted by $\partial \Omega_N$, consists of points in $\Omega_N$ which
have a neighbor not in $\Omega_N$:
\begin{equation*}
\partial \Omega_N \;=\; \big\{\, x\in \Omega_N : \exists\, y \in
(\bb Z/N)^d \setminus \Omega\;, |y-x|=1/N\, \big\}\;, 
\end{equation*}
where $\color{blue} |\,\cdot\,|$ stands for the Euclidian distance.

Let $\color{blue} \ms E$ be a subspace of $\bb R$ which represents
the possible values of the spins or occupation variables. In the case
of exclusion processes, $\ms E=\{0,1\}$. For zero-range processes,
$\ms E = \bb N \cup \{0\}$, and for KMP models, $\ms E = \bb R_+$.
The elements of $\ms E$ are denoted by the symbols $\mf x$, $\mf y$.

The state space of the microscopic dynamics is represented by
$\Sigma_N$: $\color{blue} \Sigma_N = \ms E^{\Omega_N}$ and its
elements by the Greek letters $\eta = (\eta_x : x\in \Omega_N)$,
$\xi$. Hence, $\eta_x$ stands for the value of the occupation variable
at $x\in\Omega_N$ for the configuration $\eta$.

The microscopic dynamics is composed of two pieces. The first one
describes the evolution in the bulk, while the second one the
interaction of the system with the boundary
reservoirs.

We do not discuss here the bulk dynamics nor the derivation of the
diffusivity and mobility. This has already been done in numerous
places. We refer to \cite{B8}, for example. We concentrate on the
boundary dynamics.

\subsection*{The boundary dynamics}

The system is in a mild contact with boundary reservoirs,
characterized by their chemical potentials
$\color{blue} \lambda\in\Lambda$. The boundary dynamics corresponds
to a continuous-time Markov chain taking values in $\ms E$. Its
generator, denoted by $\mc L_{\lambda}$, takes the form
\begin{equation}
\label{3-1}
(\mc L_{\lambda} f)(\mf x) \;=\; \int_{\ms E}
[f(\mf y) \,-\, f(\mf x)\,] \; r_{\lambda} (\mf x, d\mf y) \;,
\end{equation}
where $r_{\lambda} (\mf x, \,\cdot\, )$ are finite positive measures
which represent the jump rates.

For exclusion processes, $\ms E = \{0,1\}$, $\Lambda = \bb R$,
$r_{\lambda} (0, d\mf y) = [\, e^\lambda /(1+e^\lambda)\,] \,
\delta_1(d\mf y)$,
$r_{\lambda} (1, d\mf y) = [\, 1 /(1+e^\lambda)\,] \, \delta_0(d\mf
y)$, where $\delta_a(\cdot)$ stands for the Dirac measure concentrated
at $a$.  For zero-range processes, $\ms E = \{0\} \cup \bb N$,
$\Lambda = \bb R$, the jump rates are given by
$r_{\lambda} (\mf x, d\mf y) = g(\mf x) \, \delta_{\mf x -1}(d\mf y)
\,+\, \, e^\lambda \, \delta_{\mf x +1}(d\mf y)$, $\mf x\in \ms E$.
Finally, for KMP models, $\ms E = \bb R_+$, $\Lambda = (-\infty, 0)$
and
$r_{\lambda} (\mf x, d\mf y) \,=\, -\, \lambda \, e^{\lambda \mf y}\,
d\mf y$.

Of course, as boundary dynamics one could consider a Markov chain
taking values on larger spaces. For example, $\ms E^{\Delta}$ for some
finite set $\Delta$. The theory can easily be extended to this case
and this is not an important hypothesis. In Appendix \ref{sap4} we
present such a model.

Assume that for all $\lambda\in \Lambda$, the $\ms E$-valued Markov
chain induced by the generator $\mc L_\lambda$ is ergodic and has a
unique stationary state denoted by $\color{blue} m_\lambda$. Assume,
furthermore, that the measures $m_\lambda$ form an exponential family:
\begin{equation}
\label{1-01}
m_\lambda (d \mf x) \;=\; \frac{1}{Z(\lambda)} \, e^{\lambda\, \mf x
\,-\, H(\mf x)}  \, \mf n(d \mf x)
\end{equation}
for some energy $H : \ms E \to \bb R$. In this formula, $Z(\lambda)$
is the normalization constant which turns $m_\lambda$ into a
probability measure, and $\mf n$ the counting measure
($\mf n(\mf x) =1$ for all $\mf x\in\ms E$ if $\ms E$ is discrete or
$\mf n$ is the Lebesgue measure if $\ms E$ is continuous). The
function $Z(\cdot)$ is called \emph{the partition function}.

Denote by $\color{blue} c (\ms E)$ the convex envelope of $\ms E$.
Let $R: \Lambda \to c (\ms E)$ be the mean of the measure $m_\lambda$:
\begin{equation}
\label{3-2}
R(\lambda) \;=\; \int_{\ms E} \mf x\; m_\lambda (d\mf x)\;.
\end{equation}
Clearly, by definition of the partition function, $R(\lambda) =
(d/d\lambda) \log Z(\lambda)$. Taking a second derivative yields that
$R'(\lambda)$ is the variance of $\mf x$ under $m_\lambda$. In
particular, $R'(\lambda)$ is strictly positive and $R$ invertible.  
Let $\Xi\colon c(\ms E) \to \Lambda$ be the inverse of $R$:
$\color{blue} \Xi = R^{-1}$. We present in the appendix explicit
formulas for $Z$, $R$ and $\Xi$ in each model.

\subsection*{The boundary Hamiltonian}

Denote by
$\ms M^{\rm bd}_{\lambda, \rho} : c (\ms E) \times \bb R \to \bb R$,
$\lambda\in \Lambda$, $\rho\in c(\ms E)$, the function given by
\begin{equation}
\label{3-02}
\ms M^{\rm bd}_{\lambda, \rho} (p) \;=\; \int_{\ms E} \frac{1}{U_p}\,
\mc L_{\lambda} U_p  \; d m_{\Xi(\rho)} \;,
\end{equation}
where $U_p$ is the function $U_p(\mf x) = e^{p\, \mf x}$. The function
$\ms M^{\rm bd}_{\lambda, \rho} $ may take the value $+\infty$ for
certain values of $p$.

In view of formula \eqref{3-1} for the generator,
\begin{equation}
\label{x5-ab}
\ms M^{\rm bd}_{\lambda, \rho} (p) \;=\; \int_{\ms E\times \ms E}
m_{\Xi(\rho)} (d\mf x) \, r_\lambda (\mf x, d\mf y)\,
\big[\, e^{p(\mf y - \mf x)} - 1\,\big] \;.
\end{equation}
A change of variables $\mf y' = \mf y - \mf x$ yields that
\begin{equation*}
\ms M^{\rm bd}_{\lambda, \rho} (p) \;=\; \int_{\ms E}
m_{\Xi(\rho)} (d\mf x) \, \int_{\ms E^-} r_\lambda (\mf x, \mf x+ d\mf y)\,
\big[\, e^{p\mf y} - 1\,\big] \;,
\end{equation*}
provided $\ms E^- = \{\mf y - \mf x : \mf y \,,\, \mf x\in \ms E\}$.
Hence, changing the order of integrations,
\begin{equation}
\label{x5-aa}
\ms M^{\rm bd}_{\lambda, \rho} (p) \;=\;
\;=\; \int_{\ms E^-} \big[\, e^{p\mf y} - 1\,\big] \, m_{\lambda,
\rho}(d\mf y)\;,
\end{equation}
where
\begin{equation*}
m_{\lambda, \rho}(d\mf y) \;=\; \int_{\ms E}
m_{\Xi(\rho)} (d\mf x) \, r_\lambda (\mf x, \mf x+ d\mf y)\;.
\end{equation*}

\begin{example}
\label{ex1}
If the generator $\mc L_{\lambda}$ induces a Markov chain on
$\color{blue} \bb N_0 := \bb N \cup \{0\}$ or on $\{0, \dots, M\}$,
$M\ge 1$, with nearest-neighbor jumps, (it only jumps from $k$ to
$k \pm 1$), as in the case of zero-range or exclusion processes,
\begin{equation}
\label{3-03}
\ms M^{\rm bd}_{\lambda, \rho} (p) \;=\; C_\lambda(\rho) \,
[e^p-1] \;+\; A_\lambda(\rho)  \, [e^{-p}-1]  \;,
\end{equation}
where $C_\lambda(\rho)$, $A_\lambda(\rho)$ stand for the creation and
annihilation rates, respectively:
\begin{equation*}
C_\lambda(\rho) \;:=\; E_{m_{\Xi(\rho)}} [\, r_\lambda (k, k+1)\,]\;, \quad
A_\lambda(\rho) \;:=\; E_{m_{\Xi(\rho)}} [\, r_\lambda (k, k-1)\,] \;.
\end{equation*}
In this formula, $r_\lambda (k, k\pm 1)$ represents for the rate at which
the Markov chain jumps from $k$ to $k\pm 1$. The variable $k$ is
integrated with respected to the measure $m_{\Xi(\rho)}$.
\end{example}

For the simple exclusion process, $m_{\Xi(\rho)}$ is the Bernoulli
measure of parameter $\rho$, and $r_\lambda (0, 1) = R(\lambda)$,
$r_\lambda (1, 0) = 1- R(\lambda)$. Thus,
\begin{equation*}
\ms M^{\rm bd}_{\lambda, \rho} (p) \;=\; [1-\rho]\, R(\lambda) \,
[e^p-1] \;+\; \rho\, [\,1-R(\lambda)\,] \, [e^{-p}-1]  \;.
\end{equation*}
For the zero-range dynamics, $r_\lambda (k, k+1) = e^\lambda$,
$r_\lambda (k, k-1) = g(k)$. Thus,
\begin{equation*}
\ms M^{\rm bd}_{\lambda, \rho} (p) \;=\;
e^\lambda \, [e^p-1] \;+\; \Xi(\rho) \, [e^{-p}-1]  
\end{equation*}
because $E_{m_{\Xi(\rho)}} [g(k)] = \Xi(\rho)$.

The KMP model does not fall in the class described above. Here,
$m_{\Xi(\rho)}$ is the exponential measure in $\bb R_+$ with density
$\rho$, and an elementary computation yields that
\begin{equation*}
\ms M^{\rm bd}_{\lambda, \rho} (p) \;=\;
\frac{\tau}{\rho+\tau} \,\Big(\frac{1}{1-\tau \, p}\,-\, 1\Big) \;+\;
\frac{\rho}{\rho+\tau} \,\Big(\frac{1}{1+\rho \, p}\,-\, 1\Big)\;, \;\;
0<p<\tau^{-1}\;, 
\end{equation*}
where $\tau = R(\lambda)$, and
$\ms M^{\rm bd}_{\lambda, \rho} (p) \,=\, \infty$ if
$p\not \in (0, \tau^{-1})$.

\begin{remark}
\label{rm0}
The zero-range process is usually parameterized by
$\varphi = e^\lambda$. In the KMP model, keep in mind that the
chemical potential $\lambda$, used to parametrize the exponential
distributions, is negative.
\end{remark}

\begin{remark}
\label{rm1}
The reader may have recognized in \eqref{3-02} the building block of
the Donsker-Varadhan large deviations rate function:
\begin{equation*}
\ms M^{\rm bd}_{\lambda, \rho} (p) \;=\; \int_{\ms E} \frac{1}{U_p}\,
\mc L_{\lambda} U_p  \; d m_{\Xi(\rho)}  \;\ge\;
\inf_{u} \int_{\ms E} \frac{1}{u}\,
\mc L_{\lambda} u  \; d m_{\Xi(\rho)} \;=:\; -\, I_\lambda(m_{\Xi(\rho)})\;,
\end{equation*}
where the infimum is carried over all positive functions $u$ which
belong to the domain of the generator $\mc L_\lambda$. Here,
$I_\lambda$ stands for the Donsker-Varadhan rate functional of the
large deviations principle for the empirical measure of the
continuous-time Markov chain whose generator is $\mc L_\lambda$ \cite{v}.
\end{remark}

\begin{remark}
\label{rm-1}
Let $(X_t:t\ge 0)$ be the continuous-time, $\ms E$-valued Markov chain
induced by the generator $\mc L_\lambda$ introduced in
\eqref{3-1}. Assume that the process is reversible.  As we learned
from Donsker and Varadhan \cite{v}, to prove a large deviations
principle for the empirical measure
$t^{-1}\int_0^t \delta_{X_s}\, ds$, the jump rates $r_\lambda$ needed
to be tilted by a function $F:\ms E \to \bb R$, as
$r_{\lambda, F} (\mf x, d\mf y) = r_{\lambda} (\mf x, d\mf y)
e^{-[F(\mf y) - F(\mf x)]}$. The purpose of the tilting is to change the
stationary state. In fact, the equilibrium state of the Markov chain
with jump rates $r_{\lambda, F}$, denoted by $m_{\lambda, F}$, is
given by
$m_{\lambda, F} (d\mf x) = (1/Z_{\lambda, F}) \, e^{2F(\mf x)}\,
m_{\lambda} (d\mf x)$, where $Z_{\lambda, F}$ is a normalizing
constant.

The cost for the empirical measure $t^{-1}\int_0^t \delta_{X_s}\, ds$
to be close to the measure $m_{\lambda, F}$, denoted by
$I_\lambda (m_{\lambda, F})$ in Remark \ref{rm1}, is given by the
relative entropy of the perturbed dynamics with respect to the
original one \cite{v}:
\begin{equation}
\label{aa5}
I_\lambda (m_{\lambda, F}) \;=\; \lim_{t\to\infty}
\frac{1}{t}\, \bb E_{F}\Big[\,
\log \frac{d \, \bb P_{F}}
{d \, \bb P} \,\Big|_{\mc F_t}\,\Big] \;=\;
-\, \int \sqrt{\frac{dm_{\lambda, F}}{dm_{\lambda}}} \,
\mc L_\lambda \sqrt{\frac{dm_{\lambda, F}}{dm_{\lambda}}}\, dm_\lambda
\;.
\end{equation}
In this formula, $\bb P$, $\bb P_{F}$ represent the distribution of
the Markov chain with jump rates $r_{\lambda}$, $r_{\lambda, F}$,
respectively, $\bb E_{F}$ the expectation with respect to $\bb P_{F}$, and
$(d \bb P_{F}/d \bb P) \,|_{\mc F_t}$ the Radon-Nikodym derivative of
$\bb P_{F}$ with respect to $\bb P$ restricted to the $\sigma$-algebra
$\mc F_t = \sigma (X(s) : 0\le s\le t)$.

In the present context, only perturbations $F$ expressed as
$F(\mf x) = p\, \mf x$ appear. This means that only measures of the
form $e^{q\, \mf x}\, m_{\lambda} (d\mf x)$ are accessible. In other
words, only perturbations that change the chemical potential are
considered in the definition of the functional
$\ms M^{\rm bd}_{\lambda, \rho}$.  Moreover, and most importantly, in
contrast with \eqref{aa5}, $\log (d \bb P_{F}/d \bb P) \,|_{\mc F_t}$
is not integrated with respect to the stationary measure induced by
the perturbed dynamics associated the jump rated $r_{\lambda, F}$,
that is $m_{\lambda, F}$, but with respect to the stationary state
induced by the bulks dynamics.  This is a consequence of the fact that
the interaction with the boundaries is mild and dominated by the bulk
dynamics.

Hence, the functional $\ms M^{\rm bd}_{\lambda, \rho}$ has to be
understood as follows. There is a family of boundary dynamics indexed
by a chemical potential $\lambda$. The stationary state is represented
by $m_\lambda$. The system is perturbed to change its chemical
potential from $\lambda$ to $\lambda + q$. The cost of this
perturbation is not computed with respect to the new state but the one
induced by the prevalent bulk dynamics (The bulk dynamics prevails
over the boundary one because the interaction of the system with the
reservoirs is mild). We shall refer to this cost as the bulk-cost.
\end{remark}

\subsection*{Equilibrium free energy and pressure}

In the present context, the equilibrium free energy takes a simple
form. According to the postulates of statistical mechanics
\cite{fv17}, since the equilibrium states are product measure, the
pressure, denoted by $p(\cdot)$, and the free energy per unit of
volume, denoted by $f(\cdot)$ and obtained as the Legendre transform
of the pressure, are given by
\begin{equation*}
p(\lambda) \;:=\; \log Z(\lambda)\;, \quad
f (\rho) \;:=\; \sup_{\theta\in \Lambda} \big\{\theta\, \rho
\,-\, p(\theta)\, \big\} \;.
\end{equation*}
Clearly,
\begin{equation}
\label{1-06}
f (\rho) \;=\; \Xi(\rho)\, \rho \;-\; p(\Xi(\rho))\;,
\quad f' (\rho) \;=\;  \Xi(\rho) \;.
\end{equation}

Fix a reference chemical potential $\lambda$, and let
\begin{equation*}
f_\lambda (\rho) \;:=\; \sup_{\theta\in \Lambda} \big\{\theta\, \rho
\,-\, \log E_{m_\lambda} \big[\, e^{\theta \, \mf x}\, \big]\,
\big\}\;.  
\end{equation*}
By definition of $f(\,\cdot\,)$,
\begin{equation}
\label{aa6}
f_\lambda (\rho) \;=\;f (\rho) \;-\; \{\lambda\, \rho \,-\,
p(\lambda)\, \} \;\;\text{and}\;\; f''_\lambda (\rho) \;=\;
f'' (\rho)\;.
\end{equation}

It is well known that $f_\lambda (\cdot)$ is the large deviations rate
functional of the sequence $N^{-1} \sum_{1\le j\le N} \mf x_j$, where
$(\mf x_j : j\ge 1)$ are i.i.d.\! random variables distributed
according to $m_\lambda$.  The functional $f_\lambda$ is called the
\emph{equilibrium free energy}.  By \eqref{1-06} and \eqref{aa6},
\begin{equation*}
f_\lambda (\rho) \;=\; \big[\, \Xi(\rho) \,-\, \lambda\,\big]\, \rho
\;-\; \log \frac{Z(\Xi(\rho))}{Z(\lambda)}\;\;\; \text{and}\;\;\;
f'_\lambda (\rho) \;=\;  \Xi(\rho) \,-\, \lambda \;.
\end{equation*}

\subsection*{An identity}

We turn to some properties of the functional
$\ms M^{\rm bd}_{\lambda, \rho}$ needed in the next sections.  We
first claim that for all $\lambda$ and $\rho$,
\begin{equation}
\label{x17}
\ms M^{\rm bd}_{\lambda, \rho}  (\,  f'(\rho) - \lambda \,) \;=\;
0\;.
\end{equation}
By \eqref{1-06},
$\Xi(\rho) = f'(\rho)$. Hence, by \eqref{3-02} and \eqref{1-01},
\begin{equation*}
\ms M^{\rm bd}_{\lambda, \rho} ( p ) \;=\;
\frac{1}{Z(f'(\rho))}\, \int_{\ms E} e^{-p\, \mf x}\,
e^{f'(\rho)\, \mf x \,-\, H(\mf x)} \,
(\mc L_{\lambda} U_p)(\mf x)  \; \mf  n(d\mf x) 
\end{equation*}
for all $p$. Replacing the first $p$ by $f'(\rho) \,-\, \lambda$, the
right-hand side becomes
\begin{equation*}
\frac{1}{Z(f'(\rho))}\, \int_{\ms E} 
e^{\lambda \, \mf x \,-\,  H(\mf x)} \,
(\mc L_{\lambda} U_p)(\mf x)  \; \mf  n(d\mf x) \;=\;
\frac{Z(\lambda)}{Z(f'(\rho))}\, \int_{\ms E} 
(\mc L_{\lambda} U_p)(\mf x)  \; m_\lambda (d\mf x) \;.
\end{equation*}
The last term vanishes because $m_\lambda$ is the stationary state for
the dynamics induced by the generator $\mc L_{\lambda}$. This proves
\eqref{x17}.

\subsection*{The functional $\mf A_\lambda$}

Let $\mf A_\lambda$, $\lambda\in\Lambda$, be the functional
given by
\begin{equation}
\label{1-08x}
\mf A_{\lambda} (\rho,p)  \;:=\;
\kappa \,\Big\{\, \ms M^{\rm bd}_{\lambda, \rho} (p)
\;-\; \ms M^{\rm bd}_{\lambda, \rho} (0) 
\;-\; p\,  (\, \ms M^{\rm bd}_{\lambda, \rho}\,)' (0) \,\Big\}  \;.
\end{equation}
By \eqref{x5-ab}, the second term on the right-hand side,
$\ms M^{\rm bd}_{\lambda, \rho} (0)$, vanishes. It has been included
to underline that $\mf A_{\lambda} (\rho, p)$ is a first order Taylor
expansion.  By \eqref{x5-ab},
\begin{equation*}
\mf A_{\lambda} (\rho,p)  \;=\; \kappa\, 
\int_{\ms E\times \ms E} m_{\Xi(\rho)} (d\mf x) \,  r_\lambda (\mf x, d\mf
y)\, \big[\, e^{p (\mf y - \mf x)} \,-\, 1\, -\,
p\, (\mf y - \mf x)\, \big] \;\ge\; 0 \;.
\end{equation*}
In particular, in the second variable, the functional $\mf A_\lambda$
behaves quadratically close to zero:
\begin{equation*}
\mf A_\lambda (\rho,p) \;\approx \; p^2\;, \;\; p\to 0\;,
\end{equation*}
for all $\rho$ small. On the other hand, under the hypotheses of the
Example \ref{ex1}, the functional $\mf A_{\lambda}$ takes the form
\begin{equation*}
\mf A_{\lambda} (\rho,p) \;=\; C_\lambda(\rho) \,
[e^p-1-p] \;+\; A_\lambda(\rho)  \, [e^{-p}-1+p]  \;,
\end{equation*}

\subsection*{Boundary condition}

By \eqref{x5-ab}, 
\begin{equation}
\label{x16}
(\, \ms M^{\rm bd}_{\lambda, \rho} \,)' \, (0)\;=\;
\<\, 1 \,,\,  \mc L_\lambda  \, \mf x \,\>_{m_{\Xi(\rho)}} \;,
\end{equation}
where $\<\, \cdot  \,,\, \cdot \,\>_{\nu}$ represents the scalar
product in $L^2(\nu)$.

\subsection*{The reversible case}

Up to the end of this section, \emph{assume that the operator
$\mc L_\lambda$ is symmetric in} $L^2(m_\lambda)$. In other words
that, for all $\lambda \in \Lambda$, the dynamics induced by the
generator $\mc L_\lambda$ at the boundary is reversible for the
measure $m_\lambda$. This condition is fulfilled by the simple
exclusion, zero-range and KMP dynamics reviewed in the Appendices
\ref{sec2}--\ref{sap3}, but not by the exclusion process with
non-reversible boundary conditions presented in Section \ref{sap4}.

We derive below three properties of the functional
$\ms M^{\rm bd}_{\lambda, \rho}$ under this assumption. Recall the
definition of $R(\lambda)$ given in \eqref{3-2} and that, by
\eqref{1-06}, $f'(\rho) = \Xi(\rho) = R^{-1}(\rho)$.  By \eqref{x16}
and the reversibility of the measure $m_\lambda$ with respect to
$\mc L_\lambda$,
\begin{equation}
\label{aa2}
(\, \ms M^{\rm bd}_{\lambda, R(\lambda)} \,)' \, (0)\;=\;
\<\,  \mc L_\lambda  1 \,,\,  \, \mf x \,\>_{m_{\lambda}}
\;=\; 0\;.
\end{equation}
We turn to the reciprocal. We claim that
\begin{equation}
\label{aa4}
(\, \ms M^{\rm bd}_{\lambda, \rho} \,)' \, (0) \;\not =\; 0\;\; \text{if}\;\;
f'(\rho) \; \not =\; \lambda \;.
\end{equation}

Indeed, by \eqref{x16} and \eqref{1-01}, and since
$\Xi(\rho) = f'(\rho)$,
\begin{equation*}
(\, \ms M^{\rm bd}_{\lambda, \rho} \,)' \, (0) \;=\;
\frac{Z(\lambda)}{Z(\Xi(\rho))}\,
\<\, e^{(f'(\rho) - \lambda) \mf x } \,,\,
\mc L_\lambda \, \mf x \,\>_{m_\lambda} \;.
\end{equation*}
Since the measure $m_\lambda$ is reversible, 
\begin{equation*}
\begin{aligned}
(\, \ms M^{\rm bd}_{\lambda, \rho} \,)' \, (0) 
\; & =\; -\, \frac{1}{2}\,
\frac{Z(\lambda)}{Z(\Xi(\rho))}\,
\int_{\ms E\times \ms E} m_{\lambda} (d \mf x) \,
r_\lambda (\mf x, d\mf y) \,
\Big\{\, e^{(f'(\rho) - \lambda) \mf y } \,-\,
e^{(f'(\rho) - \lambda) \mf x }\, \Big\}\, (\mf y \,-\, \mf x\,) \\
& =\; -\, \frac{1}{2}\,
\int_{\ms E\times \ms E} m_{\Xi(\rho)} (d \mf x) \,
r_\lambda (\mf x, d\mf y) \,
\Big\{\, e^{(f'(\rho) - \lambda) (\mf y -\mf x)} \,-\,
1 \, \Big\}\, (\mf y \,-\, \mf x\,) \;.
\end{aligned}
\end{equation*}
Therefore, as $z\, [e^z -1] >0$ for $z\not = 0$,
\begin{equation}
\label{aa3}
[\, f'(\rho) - \lambda\,]\; 
(\, \ms M^{\rm bd}_{\lambda, \rho} \,)' \, (0)  \;<\; 0
\end{equation}
if $f'(\rho) \not = \lambda$. This proves \eqref{aa4}. \smallskip 

We conclude this section proving a last relation for
$\ms M_{\lambda, \rho}$.  We claim that
\begin{equation}
\label{x15}
\ms M_{\lambda, \rho} \big(\, f'(\rho) - f'(p)\,\big) \;=\;
\ms M_{\lambda, \rho} \big(\,  f'(p) - \lambda \,\big)
\end{equation}
for all $\lambda$, $\rho$ and $p$.

Identity \eqref{x15} asserts that the bulk-cost (in the sense of
Remark \ref{rm-1}) of changing the boundary chemical potential from
$\lambda$ to $\lambda + 2 \, [f'(\rho) - f'(p)]$ is equal to the one
of changing it from $\lambda$ to $2\, f'(p) - \lambda$. 

We turn to the proof of \eqref{x15}.  By the definition \eqref{3-02}
of $\ms M_{\lambda, \rho}$ and since $\Xi(\rho) = f'(\rho)$, the
left-hand side of \eqref{x15} is equal to
\begin{equation*}
\frac{1}{Z(f'(\rho))}\, \int_{\ms E} 
e^{f'(p)\, \mf x \,-\, H(\mf x)} \,
(\mc L_{\lambda} U_q)(\mf x)  \; \mf  n(d\mf x)
\;=\; \frac{Z(\lambda)}{Z(f'(\rho))}\, \int_{\ms E} 
U_{q'} (\mf x)  \, (\mc L_{\lambda} U_q)(\mf x)  \; m_\lambda (d\mf
x)\;, 
\end{equation*}
where $q = f'(\rho) - f'(p)$, $q' =  f'(p)- \lambda$. As the measure
$m_\lambda$ is reversible for the dynamics induced by $\mc
L_{\lambda}$, the previous expression is equal to
\begin{equation*}
\frac{Z(\lambda)}{Z(f'(\rho))}\, \int_{\ms E} 
U_{q} (\mf x)  \, (\mc L_{\lambda} U_{q'})(\mf x)  \; m_\lambda (d\mf
x) \;=\; \int_{\ms E} 
U_{q''} (\mf x)  \, (\mc L_{\lambda} U_{q'})(\mf x)  \; m_{f'(\rho)} (d\mf
x)\;,
\end{equation*}
where $q'' = \lambda - f'(p) = - \, q'$. This completes the proof of
\eqref{x15}, as $q' = f'(p)- \lambda$.

Note that \eqref{x17} follows from \eqref{x15}, but we used here the
reversibility of $m_\lambda$, while this assumption is not needed in
the derivation of \eqref{x17} presented above.

\section{The Hamiltonian formalism}
\label{sec0}

In this section, we present the thermodynamic description of
non-equilibrium driven diffusive systems in mild contact with
reservoirs. The definitions below are motivated and supported by the
microscopic dynamics reviewed in Section \ref{sec1} and in the
appendices.

Recall that $\color{blue} \Omega$ stands for the bounded domain of
$\bb R^d$ occupied by the system. The macroscopic state of the system
is described by the local density $\color{blue} \rho(x)$,
$x\in \Omega$. At each point $x$, the density $\rho(x)$ takes value in
a subset $\color{blue} \ms R$ of $\bb R$ (the set $c(\ms E)$
introduced in the previous section).  The system is in a mild contact
with boundary reservoirs, characterized by their chemical potentials
$\color{blue} \lambda\in\Lambda$, and under the action of an external
field $\color{blue} E\in \bb R^d$. The evolution is characterized by
an Hamiltonian.

\subsection*{The boundary Hamiltonian}

The boundary Hamiltonian $\ms H^{\rm bd}_{\lambda}$ is expressed in
terms of a family $\color{blue} m_{\lambda, \rho}$ of finite,
non-negative measures on $\bb R$, indexed by the chemical potential
$\lambda\in \Lambda$ and the density $\rho \in \ms R$. Let
$\ms M^{\rm bd}_{\lambda, \rho}\colon \bb R \to \overline{\bb R}$ be
the functional given by
\begin{equation}
\label{x5}
\ms M^{\rm bd}_{\lambda, \rho}  (\, p \,) \;=\;
\int_{\bb R} \big(\, e^{p\,\mf x} \,-\, 1\,\big) \; m_{\lambda,
\rho}(d\mf x)\;.
\end{equation}

The boundary Hamiltonian reads
\begin{equation}
\label{x11}
\ms H^{\rm bd}_{\lambda} \big(\, \rho\, ,\, F \, \big)  \;:=\;
\int_{\partial\Omega}
\ms M^{\rm bd}_{\lambda, \rho}  (\, F \,) \; \kappa\; d{\rm S}\;,
\end{equation}
where $\color{blue} \kappa: \bb R^d \to \bb R_+$ is a continuous,
strictly positive function wich represents the system interaction
strength with the boundary, and $d{\rm S}$ the surface measure. By
\eqref{x11},
\begin{equation}
\label{x8}
\frac{\delta \ms H^{\rm bd}_{\lambda}}{\delta F}
\big(\rho , F) (x) \;=\;
\kappa\, (\ms M^{\rm bd}_{\lambda, \rho})'  (F(x))\;, \quad
x\,\in\, \partial\Omega \;.
\end{equation}

\begin{remark}
As $\kappa$ is fixed we omit from the notation the dependence of the
boundary Hamiltonian $\ms H^{\rm bd}_{\lambda}$ on $\kappa$.
\end{remark}

\subsection*{The Hamiltonian}

The evolution of the density is described by the Hamiltonian
$\ms H_{E,\lambda}$ which takes the form
\begin{equation}
\label{3-00}
\begin{gathered}
\ms H_{E,\lambda} (\rho, F) \; =\; \ms H^{\rm bulk}_{E} (\rho, F)
\; +\; 
\ms H^{\rm bd}_{\lambda} (\rho, F )  \;, \\
\ms H^{\rm bulk}_{E} (\rho, F) \;=\; 
-\, \int_{\Omega} D(\rho) \, \nabla \rho \, \cdot \, \nabla F\; dx
\;+\; \int_{\Omega}  \sigma(\rho) \, \big\{\, E \,+\, \nabla F\, \big\} \, \cdot \, \nabla F
\; dx  \;.
\end{gathered}
\end{equation}
The \emph{diffusion coefficient} $D(\rho)$ and the \emph{mobility}
$\sigma(\rho)$ are $d\times d$ positive, symmetric matrices.  The
transport coefficients $D$ and $\sigma$ satisfy the local Einstein
relation
\begin{equation}
\label{0-1}
D(\rho) \;=\; \sigma(\rho) \, f''(\rho)
\end{equation}
where $f$ is the equilibrium free energy of the homogeneous system.
\emph{The pair $\rho$, $F$ plays the role of position and momenta},
respectively, in the Hamiltonian formalism of classical mechanics.

\begin{remark}
\label{rm2}
While the bulk Hamiltonian is expressed in terms of two
thermodynamical features, the diffusivity $D$ and the mobility
$\sigma$, the boundary Hamiltonian is written by means of a family of
measures. In all examples presented at the end of the article, $D$ and
$\sigma$ are scalars.
\end{remark}

Denote by $J_E(\rho)$ the current of the density profile $\rho$, given
by
\begin{equation}
\label{3-20}
J_E(\rho) \;=\;  -\, D(\rho) \, \nabla \rho \;+\; \sigma(\rho) \, E \;,
\end{equation}

\begin{remark}
In \cite{bdgjl14}, the bulk Hamiltonian is defined by
\begin{equation}
\label{aa1}
\ms H^{\rm bulk}_{E} (\rho, F) \;=\; 
\int_{\Omega} \big\{\,  \nabla F \,\cdot\,
\sigma(\rho) \, \nabla F \;-\; F\, 
\nabla \,\cdot\, J_E(\rho)\,\big\}\; dx\;,
\end{equation}
Of course, one could adopt this formulation, and modify accordingly
the boundary Hamiltonian to take into account the new term resulting
from the integration by parts.

However, if one adopts the definition \eqref{aa1}, the boundary
Hamiltonian will contain terms with the derivative of $\rho$. This is
not the case with the definition adopted here. The boundary
Hamiltonian, given in equation \eqref{x11}, only contains terms with
$\rho$, and not its derivative.

From a microscopic point of view, the definition \eqref{x11} is more
natural \cite{FGLN2021, BEL21}. Moreover, in the proof of the large
deviations, it is possible to show that density profiles with infinite
energy can be discarded, and one can restrict the analysis to density
profiles with generalized derivatives in $L^2$. As these profiles are
Lipschitz continuous in dimension $1$, the value of the profile at the
boundary is well defined. In contrast, it is not clear how to define
the value of the derivative of a profile at the boundary.
\end{remark}

\subsection*{Properties of the boundary Hamiltonian}

We assume that for all $\lambda$ and $\rho$,
\begin{equation}
\label{x13}
\ms M^{\rm bd}_{\lambda, \rho}  (\,  f'(\rho) - \lambda \,) \;=\;
0\;.
\end{equation}
This is property \eqref{x17} of the previous section.  We assume,
furthermore, that
\begin{equation}
\label{3-06}
\begin{gathered}
\frac{\delta \ms H^{\rm bd}_\lambda}{\delta\rho}(\rho,0) \;=\; 0 \;\;
\text{for all}\; \rho\;, \quad
[\, f'(\rho) \,-\, \lambda\,]\,
\frac{\delta \ms H^{\rm bd}_\lambda}{\delta F}(\rho,0) \;\le \; 0 \;,\\
\;\;\text{and that}\;\;
\frac{\delta \ms H^{\rm bd}_\lambda}{\delta F}(\rho,0) \;=\; 0 \;\;
\text{if and only if}\;\;  f'(\rho) \;=\; \lambda \;.
\end{gathered}
\end{equation}
The second and third conditions correspond to \eqref{aa2}, \eqref{aa4}
and \eqref{aa3} since, by \eqref{x8},
$(\delta \ms H^{\rm bd}_\lambda/\delta F) (\rho,0) = \kappa\, (\ms
M^{\rm bd}_{\lambda, \rho})' (0)$. Note that the inequality is strict
if $f'(\rho) \not = \lambda$ in view of the last property.

We turn to the assertion in \eqref{3-06} concerning the derivative
$\delta \ms H^{\rm bd}_\lambda / \delta \rho$. By the definition
\eqref{x11} of the boundary Hamiltonian, by equation \eqref{3-1} for
the generator and \eqref{1-01} for the measure,
\begin{equation*}
\frac{\delta \ms H^{\rm bd}_\lambda}{\delta \rho} (\rho, F) \;=\;
\kappa\, \Xi'(\rho)\, \Big\{\,
\big\< \, \mf x \, e^{-p\, \mf x} \,,\,
\mc L_\lambda \, e^{p\mf x} \, \big\>_{m_{\Xi(\rho)}}
\;-\; A(\rho) \,
\big\< \, \mf x \, e^{-p\, \mf x} \,,\,
\mc L_\lambda \, e^{p\mf x} \, \big\>_{m_{\Xi(\rho)}}\,\Big\}\;,
\end{equation*}
where $A(\rho) = Z'(\Xi(\rho))/Z(\Xi(\rho))$ and $p$ has to be
replaced by $F$ at the end of the computation. Hence, at $p=F=0$,
since $\ms L_\lambda 1=0$, this expression vanishes, as claimed in
\eqref{3-06}.

\subsection*{Hamilton's equation of motion}

The evolution of the pair $(\rho,F)$ is described by the Hamilton's
equation:
\begin{equation*}
\partial_t \rho \;=\; \frac{\delta \ms H_{E,\lambda}}{\delta F} \;, \quad
\partial_t F \;=\;  -\, \frac{\delta \ms H_{E,\lambda}}{\delta \rho}\;\cdot
\end{equation*}
The explicit formula for the Hamiltonian \eqref{3-00}, the divergence
theorem and the symmetry of the matrix $D$ yield the pair of equations
\begin{equation*}
\left\{
\begin{aligned}
& \partial_t u \;=\; - \, \nabla \cdot J_E(u) \;-\;
2\, \nabla \cdot \big\{ \sigma(u) \, \nabla F\, \big\} \;, \\
& \partial_t F \;=\; -\, {\bf Tr } \, \big[\, D(u) \, \text{Hess}\,
F\,\big] \,-\, \sigma' (u) \, \big[\, E \,+\, \nabla F\,
\big]\,\cdot\, \nabla F  \;.
\end{aligned}
\right.
\end{equation*}
Here, $J_E$ is the current, introduced in \eqref{3-20},
$\text{Hess } F$ stands for the Hessian matrix of $F$ and
${\bf Tr } A$ for the trace of a matrix $A$.  These equations are
complemented with the boundary conditions
\begin{equation*}
\left\{
\begin{aligned}
& \big\{ J_E (u) \;+\; 2\, \sigma(u) \, \nabla F\, \big \}
\cdot \bs n \;=\; -\, 
\frac{\delta \ms H^{\rm bd}_{\lambda}}{\delta F}(u,F)
\;, \\
& D(u) \, \nabla F \cdot \bs n \;=\; 
\frac{\delta \ms H^{\rm bd}_{\lambda}}{\delta\rho}(u,F)
\;, 
\end{aligned}
\right.
\end{equation*}
where $\bs n$ stands for the outer normal vector to $\partial \Omega$.

These equations are derived by taking the time derivative of the
equation $\ms H_{E,\lambda} (u_t, F_t) = C_0$ and integrating by
parts.  By \eqref{3-06}, the pair $(u(t), 0)$ is a solution for
Hamilton's equation of motion provided $u$ solves the hydrodynamic
equation
\begin{equation}
\label{3-08}
\left\{
\begin{aligned}
& \partial_t u \;+\; \nabla \cdot J_E (u)  \;=\; 0 \;, \\
& J_E (u) \cdot \bs n \;=\; -\, 
\frac{\delta \ms H^{\rm bd}_\lambda}{\delta F}(u,0)
\;\Big( \, = \; -\kappa\, (\ms M^{\rm bd}_{\lambda, u})'  (0) \,
\Big) \;.
\end{aligned}
\right.
\end{equation}
The last identity follows from \eqref{x8}.

\begin{remark}
Letting $\kappa \to 0$, $+\infty$ yield to Neumann boundary
conditions, $J_E (u) \cdot \bs n \,=\, 0$, and Dirichlet boundary
conditions, $(\ms M^{\rm bd}_{\lambda, u})' (0) = 0$, respectively.
\end{remark}

\begin{remark}
Equation \eqref{x16} provides an alternative formula for the current
at the boundary for the solutions of the hydrodynamic equation
\eqref{3-08}.  By \eqref{x8} and \eqref{x16},
\begin{equation*}
\frac{\delta \ms H^{\rm bd}_{\lambda}}{\delta F}
\big(\rho , 0) \;=\; \kappa\, \<\, 1 \,,\,  \mc L_\lambda  \, \mf x
\,\>_{m_{\Xi(\rho)}}  \;.
\end{equation*}
By \eqref{3-08}, the left-hand side of this equation (with a minus
sign) is equal to the value of the current at the boundary.
Therefore, the current at the boundary of the solutions of the
hydrodynamic equation can also be written as
$-\, \kappa\, \<\, 1 \,,\, \mc L_\lambda \, \mf x
\,\>_{m_{\Xi(\rho)}}$.
\end{remark}

Assume that equation \eqref{3-08} admits a unique stationary solution,
denoted by $\color{blue} \bar\rho_{E,\lambda}$. It solves the
elliptic equation
\begin{equation}
\label{3-09}
\left\{
\begin{aligned}
& \nabla \cdot \, J_E (\rho) \;=\;  0\;, \\
&  J_E (\rho) \cdot \bs n \;=\;
-\, \frac{\delta \ms H^{\rm bd}_\lambda}{\delta F}(\rho,0) \;.
\end{aligned}
\right.
\end{equation}

Assume, furthermore, that $\bar\rho_{E,\lambda}$ is an
attractor for the dynamical system induced by \eqref{3-08}.
Therefore, if $u^{(\rho)}(t)$ represents the solution of the
hydrodynamic equation \eqref{3-08} with initial condition $\rho$,
$u^{(\rho)}(0, \cdot) = \rho(\cdot)$, for every density profile
$\gamma$,
\begin{equation}
\label{3-10}
\lim_{t\to\infty} u^{(\gamma)}(t) \;=\; \bar\rho_{E,\lambda}\;.
\end{equation}

\subsection*{The action functional}

Denote by $\ms L_{E,\lambda}$ the Lagrangian associated to the
Hamiltonian $\ms H_{E,\lambda}$: For a density profile $\rho$ and a
function $G$,
\begin{equation}
\label{3-17}
\ms L_{E,\lambda} (\rho, G) \;=\; \sup_{F} \Big\{\,
\int_{\Omega} G \, F\; dx \;-\; \ms H_{E,\lambda} (\rho ,F)\,
\Big\}\;.
\end{equation}
The action functional on an interval $[T_1,T_2]$, denoted by
$I^{E,\lambda}_{[T_1,T_2]}$, is given by
\begin{equation}
\label{3-11}
I^{E,\lambda}_{[T_1,T_2]} (u) \;=\;
\int_{T_1}^{T_2} \ms L_{E,\lambda} (u(t), \partial_t u(t))
\; dt\;, 
\end{equation}
for a trajectory $u(t)$ (for each $t\ge 0$, $u(t)$ is a density
profile). The action functional indicates the cost of a path $u(t)$ in
a time interval.

\subsection*{The quasi-potential}

The quasi-potential associated to the Hamiltonian $\ms H_{E,\lambda}$,
represented by $V_{E,\lambda}$, is given by
\begin{equation}
\label{3-07} 
V_{E,\lambda} (\gamma) \;=\; \inf_{u} \;
I^{E,\lambda}_{(-\infty,0]} (u) \;,
\end{equation}
where the infimum is carried over all paths $u(t)$ starting from the
attractor $\bar\rho_{E,\lambda}$ and ending at $\gamma$:
$\lim_{t\to - \infty} u(t) = \bar\rho_{E,\lambda}$, $u(0) = \gamma $.
The quasi-potential $V_{E,\lambda}(\gamma)$ measures the minimal cost
to produce a density profile $\gamma$ starting from the stationary
profile $\bar\rho_{E,\lambda}$.

\subsection*{Hamilton-Jacobi equation}

Classical arguments in mechanics \cite{A89} imply that the
quasi-potential $V_{E,\lambda}$ solves the Hamilton-Jacobi equation
\begin{equation}
\label{3-12}
\ms H_{E,\lambda} \Big(\, \rho\, ,
\frac{\delta V_{E,\lambda}}{\delta \rho} (\rho) \,\Big) \;=\; 0\;.
\end{equation}

\subsection*{The equilibrium quasi-potential}

A state $(E,\lambda)$ is said to be an \emph{equilibrium state} if
\begin{equation}
\label{x9}
J_E(\bar\rho_{E,\lambda}) \;=\; 0\;.
\end{equation}
In this case, by \eqref{3-20} and the Einstein relation \eqref{0-1},
and by \eqref{3-09} and \eqref{3-06},
\begin{equation}
\label{bl02}
E \;=\; \nabla f'(\bar\rho_{E,\lambda}) \;\;\text{on}\;\; \Omega
\;\;\text{and}\;\;
f'(\bar\rho_{E,\lambda}) \;=\;\lambda \;\;\text{on}\;\;
\partial \Omega\;.
\end{equation}
Note that the equilibrium states in the case of strong and mild
interactions with the reservoirs are the same.

\begin{remark}
By \eqref{3-09} and \eqref{3-06}, in non-equilibrium,
$f'(\bar\rho_{E,\lambda}) \not = \lambda$ at the boundary.  This is in
sharp contrast with diffusive systems in strong interaction with
reservoirs, where $f'(\bar\rho_{E,\lambda}) = \lambda$ at the
boundary \cite{bdgjl11, bgjl1}.
\end{remark}

We claim that in equilibrium,
\begin{equation}
\label{x10}
\frac{\delta V_{E,\lambda}}{\delta \rho} (\rho)  \;=\;
f'(\rho) \;-\; f'(\bar\rho_{E,\lambda})\;,
\end{equation}
so that
\begin{equation}
\label{x12}
V_{E,\lambda} (\rho) \;=\;
\int_{\Omega} \big\{\, f(\rho) - f(\bar\rho_{E,\lambda}) -  
f'(\bar\rho_{E,\lambda}) \, \big( \rho - \bar\rho_{E,\lambda}\big) \,\big\}
\; dx \;.
\end{equation}

We turn to the derivation of \eqref{x10}. Since
$J_E(\bar\rho_{E,\lambda}) = 0$, by \eqref{bl02}, we may replace on the
right-hand side of \eqref{3-00} $E$ by
$\nabla f'(\bar\rho_{E,\lambda})$ to get that
\begin{equation*}
\ms H^{\rm bulk}_{E} \big(\, \rho\,,\, f'(\rho) - f'(\bar\rho_{E,\lambda})
\,\big) \;=\; 
-\, \int_{\Omega} D(\rho) \, \nabla \rho \, \cdot \, \nabla F\; dx
\;+\; \int_{\Omega}  \sigma(\rho) \, \nabla f'(\rho) \,  \, \cdot \, \nabla F
\; dx  \;,
\end{equation*}
where $F= f'(\rho) - f'(\bar\rho_{E,\lambda})$. By Einstein relation
\eqref{0-1}, we conclude that the right-hand side vanishes.

On the other hand, by \eqref{bl02} and \eqref{x13},
\begin{equation*}
\ms H^{\rm bd}_{\lambda} \big(\, \rho\, ,\,
f'(\rho) - f'(\bar\rho_{E,\lambda})  \, \big)  \;=\;
\ms H^{\rm bd}_{\lambda} \big(\, \rho\, ,\,
f'(\rho) - \lambda  \, \big) 
\;=\;0\;.
\end{equation*}

It follows from the two previous displayed equation that $f'(\rho) -
f'(\bar\rho_{E,\lambda})$ solves the Hamilton-Jacobi equation, proving
claim \eqref{x10}.

\subsection*{The adjoint Hamiltonian}

The Hamiltonian introduced at the beginning of this section derives
from an underlying microscopic dynamics. Denote by
$\color{blue} \bb P^*_{\rm st}$ the probability measure describing
the stationary evolution of the time-reversed process, which is still
a Markovian dynamics. We refer to $\bb P^*_{\rm st}$ as the adjoint
dynamics. The reader finds in section 1 of \cite{B2} a detailed
description of the adjoint dynamics of a Markov process.

Assume that the empirical density of the adjoint dynamics satisfies a
large deviations principle described by a Hamiltonian
$\color{blue} \ms H^\dagger_{E,\lambda}$ of the same nature as
$\ms H_{E,\lambda}$.  Denote by
$\color{blue} \ms L^\dagger_{E,\lambda}$ the Lagrangian corresponding
to the Hamiltonian $\ms H^\dagger_{E,\lambda}$.

Fix $T>0$ and a trajectory $u(t)$, $0\le t\le T$. Let $v(t) = u(-t)$.
By equation (2.2) in \cite{B2},
\begin{equation*}
V_{E,\lambda}(u(0)) \;+\; \int_0^T \ms L^\dagger_{E,\lambda} (\,
u(t)\,,\, \partial_t u (t)\,)\; dt \;=\;
V_{E,\lambda}(u(T)) \;+\; \int_{-T}^0 \ms L_{E,\lambda} (\,
v(t)\,,\, \partial_t v (t)\,)\; dt\;.
\end{equation*}
Dividing this identity by $T$ and letting $T\to 0$ yields that
\begin{equation*}
\ms L^\dagger_{E,\lambda} (\, u(0)\,,\, \partial_t u (0)\,) \;=\;
\ms L_{E,\lambda} (\, v(0)\,,\, \partial_t v (0)\,) \;+\; 
\frac{\delta V_{E,\lambda}}{\delta \rho} (u(0)) \, \partial_t u (0)\;.
\end{equation*}
Since $v (0) \,=\, u (0)$ and
$\partial_t v (0) \,=\, -\, \partial_t u (0)$, for all $\rho$,
$\gamma$,
\begin{equation}
\label{3-18}
\ms L^\dagger_{E,\lambda} (\,\rho, \gamma \,) \;=\;
\ms L_{E,\lambda} (\, \rho, -\, \gamma \,) \;+\; 
\frac{\delta V_{E,\lambda}}{\delta \rho} (\rho) \, \gamma \;.
\end{equation}

As the Hamiltonian is the convex conjugate of the Lagrangian, an
elementary computation yields that
\begin{equation}
\label{3-19}
\ms H^\dagger_{E,\lambda} (\,\rho, F \,) \;=\;
\ms H_{E,\lambda} \Big(\, \rho \,,\, \frac{\delta V_{E,\lambda}}{\delta \rho}
(\rho) \,-\, F \,\Big) \;.
\end{equation}
This formula coincides with equation (4.15) presented in
\cite{bdgjl14} for the adjoint Hamiltonian in the case of strong
boundary interactions.

The adjoint Hamiltonian $\ms H^\dagger_{E,\lambda}$ plays a central
role in the macroscopic fluctuation theory. It is shown in
\cite{BEL21} that the solution of the variational problem
\eqref{3-07}, which defines the quasi-potential, is the time-reversed
trajectory of the Hamilton's equation of motion induced by the adjoint
Hamiltonian.

Set
$\color{blue} \ms H^{\dagger, {\rm bulk}}_{E,\lambda} (\,\rho, F \,)
\,=\, \ms H^{\rm bulk}_{E} (\, \rho \,,\, \ms V(\rho) \,-\, F \,)$,
where $\ms V(\rho) = (\delta V_{E,\lambda}/\delta \rho) (\rho)$ and
define $\color{blue} \ms H^{\dagger, {\rm bd}}_{E,\lambda}$ in a
similar way. Note that $\ms H^{\dagger, \rm bulk}_{E, \lambda}$
depends on the chemical potential $\lambda$ because so does the
quasi-potential $V_{E, \lambda}$.

\subsection*{The adjoint current}

In view of the definition \eqref{3-20} of the current, the bulk
Hamiltonian $\ms H^{\rm bulk}_{E}$ can be expressed as
\begin{equation}
\label{3-24}
\ms H^{\rm bulk}_{E} (\rho, F) \; =\; 
\int_{\Omega} J_E(\rho) \, \cdot \, \nabla F\; dx
\;+\; \int_{\Omega}  \nabla F\,  \cdot \, \sigma(\rho) \, \nabla F
\; dx   \;.
\end{equation}
By \eqref{3-19}, with the bulk Hamiltonian instead of the full one, 
\begin{equation*}
\begin{aligned}
\ms H^{\dagger, {\rm bulk}}_{E,\lambda} (\,\rho, F \,) 
\; & =\; 
\int_{\Omega} J^\dagger_{E, \lambda} (\rho) \, \cdot \, \nabla F\; dx
\;+\; \int_{\Omega}  \nabla F\,  \cdot \, \sigma(\rho) \, \nabla F
\; dx  \\
&+ \; \int_{\Omega} \Big\{\, J_E(\rho) \,+\, \sigma(\rho) \, \nabla\,
\frac{\delta V_{E,\lambda}}{\delta \rho}\, (\rho)\, \Big\} \cdot 
\nabla\,
\frac{\delta V_{E,\lambda}}{\delta \rho}\, (\rho) \; dx \;,
\end{aligned}
\end{equation*}
provided we set
\begin{equation}
\label{3-22}
J^\dagger_{E, \lambda} (\rho) \;=\; -\, J_E(\rho) \;-\; 2\, \sigma (\rho)\, \nabla\, 
\frac{\delta V_{E,\lambda}}{\delta \rho}\, (\rho) \;\cdot
\end{equation}
Hence, as a function of the second variable, up to an additive
constant, the adjoint bulk Hamiltonian has the same structure as the
original one, provided we replace the current $J_E(\rho)$ by
$J^\dagger_{E, \lambda} (\rho)$.  As observed before,
$J^\dagger_{E, \lambda}$ depends on $\lambda$ because so does the
quasi-potential $V_{E, \lambda}$.

\subsection*{The adjoint hydrodynamic equation}

Computing the derivatives of the adjoint Hamiltonian yields that
$(u(t),0)$ is a solution of the adjoint Hamilton's equations provided
$u(t)$ solves the equation
\begin{equation}
\label{3-21}
\left\{
\begin{aligned}
& \partial_t u \;+\; \nabla \cdot J^\dagger_{E, \lambda} (u)  \;=\; 0 \;, \\
& J^\dagger_{E, \lambda} (u) \cdot \bs n \;=\; 
\frac{\delta \ms H^{\rm bd}_\lambda}{\delta F}
\Big(\, u, \frac{\delta V_{E,\lambda}}{\delta \rho}\, (u)\, \Big) \;,
\end{aligned}
\right.
\end{equation}
called, hereafter, the adjoint hydrodynamic equation.

\subsection*{Currents}

In view of \eqref{3-22}, it is natural to define the symmetric and
anti-symmetric currents, denoted by $J^s_{E,\lambda}(\rho)$,
$J^a_{E,\lambda}(\rho)$, respectively, as
\begin{equation}
\label{3-25}
J^s_{E,\lambda} (\rho) \;=\; -\, \sigma(\rho) \, \nabla 
\frac{\delta V_{E,\lambda}}{\delta \rho} (\rho) \;, \quad
J^a_{E,\lambda} (\rho)\;=\;
(1/2) \, \big\{\, J_E(\rho) \,-\, J^\dagger_{E,\lambda}(\rho)
\,\big\}\;.
\end{equation}

The Hamilton-Jacobi equation provides an orthogonality relation between
the symmetric and anti-symmetric currents.  In equation \eqref{3-24},
replace $F$ by $\delta V_{E,\lambda}/\delta \rho$ and recall the
definition of the symmetric current to express the second term on the
right-hand side of \eqref{3-24} as a function of
$J^s_{E,\lambda}(\rho)$. As
$J_{E}(\rho) = J^s_{E,\lambda}(\rho) + J^a_{E,\lambda}(\rho)$,
\begin{equation*}
\ms H^{\rm bulk}_{E} \big(\, \rho \,,\,
\frac{\delta V_{E,\lambda}}{\delta \rho} (\rho) \,\big) \;=\;
-\, \int_{\Omega} J^a_{E,\lambda}(\rho) \cdot 
\frac{1}{\sigma(\rho)} \, J^s_{E,\lambda}(\rho) \; dx \;.
\end{equation*}
In particular, the Hamilton-Jacobi equation \eqref{3-12} becomes an
orthogonality relation between the anti-symmetric and the symmetric
currents:
\begin{equation}
\label{3-15}
\int_{\Omega} J^a_{E,\lambda}(\rho) \cdot 
\frac{1}{\sigma(\rho)} \, J^s_{E,\lambda}(\rho) \; dx
\;=\; 
\ms H^{\rm bd}_{\lambda} \big(\, \rho\, ,\,
\frac{\delta V_{E,\lambda}}{\delta \rho} (\rho)\, \big)
\;= \int_{\partial\Omega} \ms M^{\rm bd}_{\lambda, \rho}
\big(\, 
\frac{\delta V_{E,\lambda}}{\delta \rho} (\rho)\, \big) \;
\kappa \; dS \;.
\end{equation}

\begin{remark}
In the case of strong boundary interactions, it has been shown that
the symmetric and the asymmetric currents satisfy an orthogonality
relation (see equation (4.5) in \cite{bgjl1} and equation (2.22) in
\cite{bdgjl14}). This orthogonality relation is fundamental in the
proof of the Clausius inequality for the renormalized work in
\cite{bdgjl14}.

It would be interesting to obtain a similar relation in the present
context.  In order to achieve this, one would need to write the
right-hand side of \eqref{3-15} as the scalar product of the symmetric
boundary current with the anti-symmetric one. This would permit to
interpret equation \eqref{3-15} as an orthogonality relation between
the symmetric and the anti-symmetric current. We were not able to.
\end{remark}

\section{A formula for the quasi-potential}
\label{sec6}

In few cases, it is possible to derive an explicit formula for the
quasi-potential.  Assume that $d=1$ and that there is no external
field, $E=0$. Then, for each density $\rho$,
\begin{equation}
\label{3-14}
\frac{\delta V_{E,\lambda}}{\delta \rho} \;=\;
f'(\rho) \;-\; f'(F) 
\end{equation}
is a solution of the Hamilton-Jacobi equation provided $F$ solves the
equation
\begin{equation}
\label{3-13}
\left\{
\begin{aligned}
& \Delta f'(F) \;+\; \frac{\sigma(\rho) - \sigma(F)}{d(\rho) - d(F)}
\, |\, \nabla f'(F)\, |^2 \;=\; 0 \;, \\
& [\, d(\rho) - d(F)\,]\, \nabla f'(F) \cdot \bs n \;=\;
\kappa\, \ms M^{\rm bd}_{\lambda, \rho}
\big(\, f'(\rho) \;-\; f'(F) \,\big) \;.
\end{aligned}
\right.
\end{equation}
In this formula, $d$ is the primitive of $D$: $d'(\rho) = D(\rho)$.

The proof of this claim is similar to the one presented in
\cite{B2}. Let $\Gamma = f'(\rho) \;-\; f'(F)$, recall the definition
of the bulk Hamiltonian and the Einstein relation \eqref{0-1} to get
that
\begin{equation*}
\begin{aligned}
\ms H^{\rm bulk}_{E} \big(\, \rho \,,\, \Gamma  \,\big)  \; & =\; 
\int_{\Omega} \sigma(\rho) \, \nabla \Gamma \cdot \nabla \Gamma   \;dx
\; -\; \int_{\Omega} \sigma(\rho) \nabla f'(\rho) \, \cdot \, \nabla \Gamma  \;dx \\
& =\; -\, \int_{\Omega} \sigma(\rho) \, \nabla f'(\rho) \cdot \nabla f'(F)  \;dx
\;+\; \int_{\Omega} \sigma(\rho) \, |\, \nabla f'(F) \,|^2  \;dx \;.
\end{aligned}
\end{equation*}
By Einstein relation \eqref{0-1}, and since $d$ is a primitive of $D$,
$\sigma(\rho) \, \nabla f'(\rho) = D(\rho) \, \nabla \rho = \, \nabla
d(\rho)$. Therefore, the first term on the right-hand side can be
written as
\begin{equation*}
-\, \int_{\Omega} \big[\, \nabla d(\rho) \,-\, \nabla d(F) \,\big] \cdot
\nabla f'(F)  \;dx \;-\;
\int_{\Omega} \nabla d(F)  \cdot \nabla f'(F)  \;dx \;.
\end{equation*}
By the divergence theorem and since $\nabla d(F)  = \sigma(F) \,
\nabla f'(F)$, this expression is equal to
\begin{equation*}
\int_{\Omega} \big[\, d(\rho) \,-\, d(F) \,\big] \, \Delta f'(F)  \;dx
\;-\; \int_{\partial \Omega} \big[\, d(\rho) \,-\, d(F) \,\big]\, 
\nabla f'(F) \cdot \bs n \; dS
\;-\;
\int_{\Omega} \sigma(F) \, |\, \nabla f'(F) \,|^2 \;dx \;.
\end{equation*}
By \eqref{x11} and \eqref{3-13}, the expression appearing in second
integral is equal to
$\ms H^{\rm bd}_\lambda (\, \rho \,,\, f'(\rho) \;-\; f'(F))$

Up to this point, we proved that
\begin{equation*}
\begin{aligned}
\ms H_{E,\lambda} \big(\, \rho \,,\, \Gamma  \,\big)  \; 
=\; \int_{\Omega} \big[\, d(\rho) \,-\, d(F) \,\big] \, \Delta f'(F)  \;dx
\;+\; \int_{\Omega} \big[\, \sigma(\rho) \,-\, \sigma(F)\,\big]
\, |\, \nabla f'(F) \,|^2 \;dx \;.
\end{aligned}
\end{equation*}
On the left-hand side, the bulk Hamiltonian appearing at the
beginning of the computation has been replaced by the full one in view
of the last observation of the previous paragraph. The right-hand side
vanishes in view of the first equation in \eqref{3-13}. This completes
the proof of the claim. \smallskip

It has been proved in \cite{BEL21} that for symmetric exclusion
processes in mild contact with reservoirs the quasi-potential is given
by \eqref{3-14}. For zero-range models, this is also easy to check
since in this case even in nonequilibrium the stationary states are
product measures. Actually, when $\sigma$ is constant (as in
Ginzburg-Landau dynamics) or $\sigma = d$ (as in zero-range dynamics),
the first equation in \eqref{3-13} becomes autonomous in $F$.

\begin{remark}
\label{rm-x}
For the one-dimensional exclusion process and the zero-range dynamics,
the boundary conditions of equation \eqref{3-13} coincide with the
ones of equation \eqref{3-09}.  That is, the boundary conditions for
the stationary density profile in equation \eqref{3-09} and the ones
for the auxiliary function $F$ in equation \eqref{3-13} coincide. This
is the also case for interacting particle systems in strong
interaction with the boundary reservoirs, where the boundary
conditions are of Dirichlet type \cite{DLS, B2}.

To derive the identity of the boundary conditions, observe that by the
Einstein relation and since $J_E(F) = -\, D(F) \nabla F$, the boundary
condition in \eqref{3-13} can be restated as
\begin{equation*}
J(F) \cdot \bs n \;=\; - \, \kappa \,
\frac{\sigma(F)}{d(\rho) - d(F)} \,
\, \ms M^{\rm bd}_{\lambda, \rho}
\big(\, f'(\rho) \;-\; f'(F) \,\big)\;.
\end{equation*}
For zero-range and exclusion dynamics, a computation, presented at the
appendix, yields that for all $\mu$, $\varrho$ and $p$
\begin{equation}
\label{x18}
\frac{\sigma(p)}{d(\varrho) - d(p)} \,
\, \ms M^{\rm bd}_{\mu, \varrho}
\big(\, f'(\varrho) \;-\; f'(p) \,\big) \;=\;
\big(\, \ms M^{\rm bd}_{\mu, p}\,\big)' (0)\;.
\end{equation}
This proves that the boundary condition in \eqref{3-13} and
\eqref{3-09} coincide in view of \eqref{x8}.  Therefore, in these
examples equation \eqref{3-13} becomes
\begin{equation}
\label{3-131}
\left\{
\begin{aligned}
& \Delta f'(F) \;+\; \frac{\sigma(\rho) - \sigma(F)}{d(\rho) - d(F)}
\, |\, \nabla f'(F)\, |^2 \;=\; 0 \;, \\
& J(F) \cdot \bs n \;=\; -\,
\kappa\, (\ms M^{\rm bd}_{\lambda, F})' (0) \;.
\end{aligned}
\right.
\end{equation}

Note that the left-hand side of \eqref{x18} depends on $\varrho$,
while the right-hand side does not. 
\end{remark}

\begin{remark}
\label{rm-x2}
Contrarily, for the KMP model, the boundary conditions of equation
\eqref{3-13} are different from the ones of equation \eqref{3-09}.
The equation for the auxiliary function is presented in
\eqref{3-13d}. This is in sharp contrast with the case of strong
interaction with the boundary, where the boundary conditions coincide.
\end{remark}

\begin{remark}
Equation \eqref{x15} provides an identity for the term
$\ms M^{\rm bd}_{\mu, \varrho} (\, f'(\varrho) \;-\; f'(p) \,)$ which
holds under the general hypotheses of Section \ref{sec1}.
\end{remark}

\begin{remark}[Boundary conditions for the quasi-potential]
When the system interacts strongly with the boundary reservoirs, the
boundary equations for the hydrodynamic equation and for the adjoint
hydrodynamic equation are of Dirichlet type. For this reason, one may
restrict the investigation of the quasi-potential to density profiles
which satisfy Dirichlet boundary conditions. Since the boundary
conditions for the auxiliary function $F$, introduced in \eqref{3-14},
are also of Dirichlet type, in view of \eqref{3-14}, one concludes
that the functional derivative of the quasi-potential vanishes at the
boundary when computed at a density profile.

For systems in mild interaction with the boundary, the situation is
completely different. The boundary conditions for the hydrodynamic and
the adjoint hydrodynamic equations, equations \eqref{3-08} and
\eqref{3-21}, respectively, are different. There are no reasons to
restrict the attention to particular density profiles. But even
considering profiles which satisfy the Robin boundary conditions
appearing in \eqref{3-08} and \eqref{3-09}, we were not able to derive
boundary conditions for the functional derivative of the
quasi-potential computed at fixed density profiles. We did not push
this investigation too much as we did not need an equation for the
boundary condition of the functional derivative of the quasi-potential
in our analysis.
\end{remark}

In terms of the auxiliary function $F$, the currents take the form
\begin{gather*}
J^s_{0,\lambda}(\rho)  \;=\; J(\rho) \;+\; \sigma(\rho)\, \nabla
f'(F)\;, \quad
J^a_{0,\lambda}(\rho)  \;=\; -\,  \sigma(\rho)\, \nabla
f'(F)\;, \\
 J^\dagger_{0,\lambda}(\rho) \;=\; J(\rho) \;+\; 2\, \sigma(\rho)\, \nabla
f'(F)\;.
\end{gather*}

\subsection*{Orthogonal relation for the currents}

Assume that the matrix $D(\rho)$ is a scalar and that the external
field $E$ vanishes. We keep $E$ in the notation though it vanishes.
Recall the representation \eqref{3-14} for the quasi-potential. By the
Einstein relation \eqref{0-1},
\begin{equation*}
J^s_{E,\lambda}(\rho) \;=\; -\, D(\rho)\, \nabla \rho \;+\;
\sigma(\rho)\, \nabla f'(F)\; =\; J_E(\rho) \;+\;
\sigma(\rho)\, \nabla f'(F)\;.
\end{equation*}
Therefore,
$J^a_{E,\lambda}(\rho) \,=\, -\, \sigma(\rho)\, \nabla f'(F)$.  Keep
in mind that $E=0$ although it remains in the notation.  By this
formula for the anti-symmetric current and the second equation in
\eqref{3-13} the orthogonality relation \eqref{3-15} can be written as
\begin{equation}
\label{3-16}
\int_{\Omega} J^a_{E,\lambda}(\rho) \cdot 
\frac{1}{\sigma(\rho)} \, J^s_{E,\lambda}(\rho) \; dx
\;+\;  \int_{\partial\Omega} 
\frac{d(\rho) - d(F)}{\sigma(\rho)} \,
J^a_{E,\lambda}(\rho)  \cdot \bs n  \; dS \;=\;  0 \;.
\end{equation}

\section{A Clausius inequality}
\label{sec3}

Fix a time-dependent chemical potential $\lambda(t,x)$ and external
field $E(t,x)$. For a density profile $\rho$, let $u(t,x)$ be the
solution of
\begin{equation}
\label{3-08b}
\left\{
\begin{aligned}
& \partial_t u \;+\; \nabla \cdot J_{E(t)} (u)  \;=\; 0 \;, \\
& J_{E(t)} (u) \cdot \bs n \;=\; -\,
\frac{\delta \ms H^{\rm bd}_{\lambda(t)}}{\delta F}(u,0) \;, \\
& u(0,\cdot) \;=\; \rho(\cdot)\;, 
\end{aligned}
\right.
\end{equation}
where $J_{E(t)} (u)$ is given by \eqref{3-20} with $u(t)$, $E(t)$
replacing $\rho$, $E$, respectively.

The energy exchanged between the system and the external
reservoirs and fields in the time interval $[0,T]$ is given by
\begin{equation*}
\begin{aligned}
& W_{[0,T]}({\lambda (\cdot), E (\cdot), \rho})  \\
& \quad \;:=\; \int_{0}^{T} \Big\{ 
- \int_{\partial\Omega} \lambda (t,x)
\: j(t,x) \cdot \bs n(x) \; d{\rm S}(x)
\: +\int_\Omega  j(t,x) \cdot E(t,x)\; dx \Big\} \; dt\;, 
\end{aligned}
\end{equation*}
where $\color{blue} j(t,x) = J_{E(t)} (u(t)) \, (x)$ is the current
of profile $u(t, \cdot)$.  The first term on the right-hand side is
the energy provided by the reservoirs while the second is the energy
provided by the external field. We claim that
\begin{equation}
\label{03}
W_{[0,T]}({\lambda (\cdot), E (\cdot), \rho})  \;\ge\;  F (u(T)) \;-\;  F(\rho) \;,
\end{equation} 
where $F$ is the equilibrium free energy functional defined by
\begin{equation}
\label{10}
F(\rho) \;:=\;  \int_\Omega f (\rho(x))\; dx\;.
\end{equation}

Indeed, dropping from the notation the dependence on $x$, adding and
subtracting $f'(u(t))$ in the boundary term, rewrite the energy
exchanged as
\begin{equation*}
\begin{split}
W_{[0,T]}({\lambda (\cdot), E (\cdot), \rho})  \; =\; & -\, 
\int_{0}^{T}\! dt \, \int_{\partial\Omega}
[\lambda (t) - f'(u(t))] \; j(t) \cdot  \bs n \; dS \\
& +\, \int_{0}^{T}\! dt \, \Big\{ 
- \int_{\partial\Omega} f'(u(t)) \; j(t) \cdot  \bs n \; dS
\;+\; \int_\Omega  j(t) \cdot E(t) \;dx\, \Big\}\;.
\end{split}
\end{equation*}
By the divergence theorem, the right-hand side is equal to
\begin{equation}
\label{x7}
\begin{aligned}
& -\, \int_{0}^{T}\! dt \, \int_{\partial\Omega}
[\lambda (t) - f'(u(t))] \, 
\, j(t) \cdot  \bs n \; dS  \\
&\quad +\;
\int_{0}^{T} \!dt \int_\Omega  \big\{ 
- \nabla\cdot \big[ f'(u(t) ) \, j(t) \big] +  j (t) \cdot E
(t)\big\} \; dx \;.
\end{aligned}
\end{equation}

Recall from \eqref{1-08x} the definition of the functional
$\mf A_{\lambda} (\rho,\,\cdot\,)$, copied here for the reader's 
convenience:
\begin{equation}
\label{1-08xb}
\mf A_{\lambda} (\rho,p)  \;=\;
\kappa \,\Big\{\, \ms M^{\rm bd}_{\lambda, \rho} (p)
\;-\; \ms M^{\rm bd}_{\lambda, \rho} (0) 
\;-\; p\,  (\, \ms M^{\rm bd}_{\lambda, \rho}\,)' (0) \,\Big\}  \;.
\end{equation}
It has been shown in Section \ref{sec1}
that $\mf A_{\lambda}$ is positive and that close to zero it behaves quadratically in
the second variable:
\begin{equation}
\label{x2b}
\mf A_{\lambda} (\rho,p)  \;\ge\; 0\;, \quad
\mf A_{\lambda} (\rho,p) \;\approx \;  p^2 \;, \;\; p\,\to\, 0\;,
\end{equation}
for all $\rho$. 

By the boundary conditions in \eqref{3-08b} and \eqref{x8},
\begin{equation*}
-\,  [\, \lambda (t) \,-\, f'(u(t))\, ] \, 
\, j(t) \cdot  \bs n \;=\;
\kappa\, [\, \lambda (t) \,-\, f'(u(t)) \,]\, 
(\, \ms M^{\rm bd}_{\lambda (t), u(t)}\,)' (0)\;.
\end{equation*}
Since $\ms M^{\rm bd}_{\lambda, \rho} (0) = 0$ and, by \eqref{x13},
$\ms M^{\rm bd}_{\lambda, \rho} (\, f'(\rho) \,-\, \lambda\,) =0$
this expression is equal to
\begin{equation*}
\kappa \,\Big\{\, \ms M^{\rm bd}_{\lambda (t), u(t)} (p)
\;-\; \ms M^{\rm bd}_{\lambda (t), u(t)} (0) 
\;-\; p\,  (\, \ms M^{\rm bd}_{\lambda (t), u(t)}\,)' (0) \,\Big\}
\;=\; \mf A_{\lambda (t)} (u(t),p) \;,
\end{equation*}
for $p = f'(u(t)) \,-\, \lambda (t)$.  Therefore, the first term in
\eqref{x7} can be written as
\begin{equation*}
\int_{0}^{T}\! dt \, \int_{\partial\Omega}
\mf A_{\lambda (t)} ( \, u(t) \,,\, f'(u(t)) - \lambda (t) \,)
\; dS \;.
\end{equation*}
By \eqref{x2b}, this expression is positive.  On the other hand, the
second term of \eqref{x7} is equal to
\begin{equation*}
\begin{split}
&\int_{0}^{T} \!dt \int_\Omega
\big[\, - f'(u(t)) \, \nabla \cdot j(t)  \,-\,  f''(u(t)) \, \nabla u(t)
\cdot j(t) \,+\,  j(t)\cdot E(t) \, \big] \; dx   \\
&\quad = \; \int_{0}^{T} \!dt  \, \frac{d}{dt}  \int_\Omega 
f( u(t) ) \; dx
\;+\; \int_{0}^{T} \!dt  \int_\Omega 
j(t)\cdot \sigma(u(t) )^{-1} j(t) \; dx\; ,
\end{split}
\end{equation*}
where we used the Einstein relation \eqref{0-1}, and the
definition \eqref{3-20} of the current $j(t) = J_{E(t)}(u(t))$.  In
conclusion, we proved that
\begin{equation}
\label{04}
\begin{split}
W_{[0,T]}({\lambda (\cdot), E (\cdot), \rho})  \; & =\;  F (u(T)) - F(\rho)
\, +\, \int_{0}^{T} \!dt  \int_\Omega 
j(t)\cdot \sigma(u(t) )^{-1} j(t) \; dx \\
& + \; \int_{0}^{T}\! dt \, \int_{\partial\Omega}
\mf A_{\lambda (t)} ( \, u(t) \,,\,
f'(u(t)) - \lambda (t) \,)\; dS  \;.
\end{split}
\end{equation}
Since the last two term are positive, inequality \eqref{03}
follows.

This argument provides a dynamic derivation of the second law of
thermodynamics as expressed by the Clausius inequality \eqref{03}. The
key ingredients have been the assumption of local equilibrium together
with the local Einstein relationship \eqref{0-1}.

\section{Transformation along equilibrium states}
\label{sec4}

In this section, we examine transformations along equilibrium
states. Recall from \eqref{x9} that equilibrium states are
characterized by the absence of current at stationarity:
$J(\bar\rho_{E,\lambda})=0$.

\subsection*{Reversible and quasi static transformations}

Fix $T>0$, and two equilibrium states $(E_0, \lambda_0)$, $(E_1,
\lambda_1)$ so that
\begin{equation}
\label{bl02-bis}
J(\bar\rho_{E_0, \lambda_0}) \;=\;
J(\bar\rho_{E_1, \lambda_1}) \;=\; 0\;.
\end{equation}
Consider a system initially in the state
$\bar\rho_0 = \bar\rho_{E_0, \lambda_0}$ which is driven to a new
state $\bar\rho_1 = \bar\rho_{E_1, \lambda_1}$ by changing the
chemical potential and the external field in time in a way that
$(E(t), \lambda(t)) =(E_0, \lambda_0)$ for $t\le 0$ and
$(E(t), \lambda(t)) =(E_1, \lambda_1)$ for $t\ge T$.  This
transformation from $\bar\rho_0$ to $\bar\rho_1$ is called
\emph{reversible} if the energy exchanged with the reservoirs is
minimal.  A basic thermodynamic principle asserts that reversible
transformation are accomplished by a sequence of equilibrium states
and are well approximated by \emph{quasistatic} transformations,
transformations in which the variation of the thermodynamical
variables is very slow so that the density profile at time $u(t)$ is
very close to the stationary profile $\bar\rho_{E(t), \lambda(t)}$.

Let $u(t,x)$, $t \ge 0$, $x\in\Omega$, be the solution of
\eqref{3-08b} with initial condition $\rho = \bar\rho_0$.  Recall that
we denote by $j(t)$ the current at time $t$ of the density profile
$u(t)$: $j(t) = J_{E(t)}(u(t))$.  Since the thermodynamical variables
are equal to $(E_1, \lambda_1)$ for $t\ge T$, as $t\to\infty$, $u(t)$
and the current $j(t)$ relax exponentially fast to $\bar\rho_1$ and
$J_{E_1}(\bar\rho_1)$, respectively.  By \eqref{bl02-bis},
$J_{E_1}(\bar\rho_1)=0$, and, by \eqref{bl02} and \eqref{1-08xb},
$\mf A_{\lambda_1} ( \, \bar\rho_1 \,,\, f'(\bar\rho_1) -
\lambda_1 \,) = \mf A_{\lambda_1} ( \, \bar\rho_1 \,,\, 0 \,)
=0$. Therefore, the integrals in \eqref{04} are finite as $T\to\infty$
and
\begin{equation}
\label{11}
\begin{aligned}
W(\lambda (\cdot), E (\cdot), \bar\rho_0)
\; & = \;  F (\bar\rho_1) \;-\; F(\bar\rho_0)
\;+\; \int_{0}^{\infty} \!dt  \int_\Omega 
j(t)\cdot \sigma(u(t) )^{-1} j(t) \; dx \\
& + \; \int_{0}^{\infty}\! dt \, \int_{\partial\Omega}
\mf A_{\lambda (t)} ( \, u(t) \,,\, f'(u(t)) - \lambda (t) \,) 
\; dS \\
& \ge \; F (\bar\rho_1) \;-\; F(\bar\rho_0) \; . 
\end{aligned}
\end{equation}
Last inequality follows from \eqref{x2b}.  Note that we did not
assume any regularity of the thermodynamical variables in time so that
they can also be discontinuous.

It remains to show that in the quasistatic limit equality is achieved
in \eqref{11}.  That is the thermodynamic relation
\begin{equation}
\label{12}
W =  \Delta F
\end{equation}
holds, where $\Delta F = F (\bar\rho_1) \,-\, F(\bar\rho_0)$ is the
variation of the free energy.  If this is the case, by running the
transformation backward in time, we can return to the original state
exchanging the energy $-\, \Delta F$.  For this reason the
transformations for which \eqref{11} becomes an equality are called
reversible.

However, for any fixed transformation the inequality in \eqref{11} is
strict because the last two terms on the right-hand side of the
identity in \eqref{11} are strictly positive. The second one is
strictly positive in view of the last assertion of
\eqref{3-06}. Hence, reversible transformations cannot be achieved
exactly.  We can however exhibit a sequence of transformations for
which these strictly positive terms can be made arbitrarily small.
This sequence of transformations is what is called quasistatic
transformations.

Fix smooth functions $\lambda(t)$, $E(t)$ such that
$(\lambda(0), E(0)) = (\lambda_0, E_0)$,
$(\lambda(t), E(t)) = (\lambda_1, E_1)$ for $t\ge T$.  Assume that
$(\lambda(t), E(t))$ are equilibrium states for all $t\ge 0$.  This
means that $J(\bar\rho_{\lambda (t), E (t)}) =0$ for all $t\ge 0$.
Given $\delta >0$, set $\lambda_\delta (t) = \lambda (\delta t)$,
$E_\delta (t) = E (\delta t)$.  The sum of the last two terms on the
right-hand side of \eqref{11} is given by
\begin{equation*}
\begin{aligned}
& \int_{0}^{\infty} \!dt  \int_\Omega
\{\, \nabla f' (u_\delta(t))  \,-\, E_\delta(t)\, \} \cdot 
\sigma(u_\delta(t)) \, \{\, \nabla f' (u_\delta(t))
\,-\, E_\delta(t)\,\} \; dx \; \\
&\quad + \; \int_{0}^{\infty}\! dt \, \int_{\partial\Omega}
\mf A_{\lambda_\delta (t)} ( \, u_\delta(t) \,,\,
f'(u_\delta(t)) - \lambda_\delta (t) \,)\; 
dS \;,
\end{aligned}
\end{equation*}
where $u_\delta$ is the solution to \eqref{3-08b} with initial
condition $\bar\rho_0$ and parameters $\lambda_\delta(t)$,
$E_\delta(t)$.

Recall that $\bar\rho_{\lambda_\delta(t), E_\delta(t)}$ is the
equilibrium state associated to the thermodynamical variables
$\lambda_\delta (t)$, $E_\delta(t)$.  Since
$J(\bar\rho_{\lambda (t), E (t)}) =0$,
$E_\delta(t) = \nabla f' (\bar\rho_{\lambda_\delta(t), E_\delta(t)})$.
Therefore, the previous expression is equal to
\begin{equation*}
\begin{aligned}
& \int_{0}^{\infty} \!dt  \int_\Omega
\{\, \nabla f' (u_\delta(t))  \,-\,
\nabla f' (\bar\rho_{\lambda_\delta(t), E_\delta(t)}) \, \} \cdot 
\sigma(u_\delta(t)) \, \{\, \nabla f' (u_\delta(t))
\,-\, \nabla f' (\bar\rho_{\lambda_\delta(t), E_\delta(t)}) \,\} \; dx \; \\
&\quad + \; \int_{0}^{\infty}\! dt \, \int_{\partial\Omega}
\mf A_{\lambda_\delta (t)} ( \, u_\delta(t) \,,\,
f'(u_\delta(t)) -  \lambda_\delta(t) \,)\; 
dS \;,
\end{aligned}
\end{equation*}
The difference between the solution of the hydrodynamic equation
$u_\delta(t)$ and the stationary profile
$\bar\rho_{\lambda_\delta(t), E_\delta(t)}$ is of order $\delta$
uniformly in time, and so is the differences
$f' (u_\delta(t)) - f'(\bar\rho_{\lambda_\delta(t), E_\delta(t)})$.
As the integration over time essentially extends over an interval of
length $\delta^{-1}$, the first term of the previous expression
vanishes for $\delta\to 0$.  Similarly, by \eqref{x2b},
$\mf A_{\lambda_\delta (t)} ( \, u_\delta(t) \,,\,
f'(u_\delta(t)) - \lambda_\delta(t) \,)$ is bounded by
$C_0 [\, f'(u_\delta(t)) - \lambda_\delta(t) \,]^2 = C_0 [\,
f'(u_\delta(t)) - f'(\bar\rho_{\lambda_\delta(t), E_\delta(t)})
\,]^2$. Hence, the second term of the previous expression also
vanishes as $\delta\to 0$.

This implies that equality in \eqref{11} is achieved in the limit
$\delta\to 0$.  Note that in the previous argument we did not use any
special property of the path $\lambda(t)$ besides its smoothness in
time. Otherwise, the trajectory $(E(t), \lambda (t))$ from
$(E_0, \lambda_0)$ to $(E_1, \lambda_1)$ can be arbitrary along
equilibrium states.

\subsection*{Excess work}

Consider a transformation $(E(t), \lambda (t))$, $t\ge 0$, and an
initial density profile $\rho$.  Assume that
$(E(t), \lambda (t))\to (E_1, \lambda_1)$, as $t\to +\infty$ fast
enough, where $(E_1, \lambda_1)$ defines an equilibrium state (that is
$J(\bar\rho_{E_1, \lambda_1})=0$).  The \emph{excess work}
$W_\mathrm{ex} = W_\mathrm{ex}(\lambda (\cdot),E (\cdot),\rho)$ is
defined as the difference between the energy exchanged between the
system and the external driving and the work involved in a reversible
transformation from $\rho$ to $\bar\rho_1$, namely
\begin{equation}
\label{exwork} 
\begin{aligned}
W_\mathrm{ex} \;=\; W  \;-\;  \min W
\; &=\;  \int_0^\infty \!dt \int_{\Omega}  j(t) \cdot
\sigma(u(t))^{-1} j(t) \; dx \\
& + \; \int_{0}^{\infty}\! dt \, \int_{\partial\Omega}
\mf A_{\lambda (t)} ( \, u (t) \,,\,
f'(u (t)) - \lambda (t) \,)\;  dS \;,
\end{aligned}
\end{equation}
where we used \eqref{11} as well as the fact that the minimum of $W$
is given by the right-hand side of \eqref{12}.  Observe that
$W_\mathrm{ex}$ is a positive functional of the transformation
$(E(t), \lambda(t))$ and the initial condition $\rho$. Of course, by
taking a sequence of quasi-static transformations $W_\mathrm{ex}$ can
be made arbitrarily small.

\subsection*{Relaxation path  and availability}

Consider an equilibrium system in the state $\bar\rho_0$,
characterized by a chemical potential $\lambda_0$ and an external
field $E_0$. This system is put in contact with reservoirs at constant
chemical potential $\lambda_1$ and an external field $E_1$, different
from the chemical potential $\lambda_0$ and the external field $E_0$
associated to $\bar\rho_0$. Assume that $(E_1, \lambda_1)$ is an
equilibrium state.  For $t> 0$ the system thus evolves according to
the hydrodynamic equation \eqref{3-08b} with initial condition
$\bar\rho_0$, external field $E_1$, and boundary condition
$\lambda_1$.

When $t\to +\infty$ the system relaxes to the equilibrium state
$\bar\rho_1$. As the path $(E(\cdot), \lambda(\cdot))$ is constant in
time, the excess work is a function of $(E_1, \lambda_1)$ and
$\bar\rho_0$ and we denote $W_\mathrm{ex} (E(\cdot), \lambda(\cdot),
\bar\rho_0)$ simply by $W_\mathrm{ex} (E_1, \lambda_1, \bar\rho_0)$. 

In view of \eqref{exwork}, the constitutive equation \eqref{3-20} and
\eqref{11}, the excess work along such a path is given by
\begin{equation*}
\begin{aligned}
W_\mathrm{ex} (E_1, \lambda_1, \bar\rho_0) 
= &- \; \int_0^{\infty}\!dt \int_\Omega 
\big[ \nabla f'(u(t)) - E_1 \big] \cdot j(t) \; dx \\
& + \; \int_{0}^\infty\! dt \, \int_{\partial\Omega}
\mf A_{\lambda_1} ( \, u (t) \,,\,
f'(u (t)) - \lambda_1 \,)\; dS\;.
\end{aligned}
\end{equation*}
Since $(E_1,\lambda_1)$ is an equilibrium state, $J(\bar\rho_1)=0$, so
that $\nabla f'(\bar\rho_1) = E_1$. We may thus replace $E_1$ by
$\nabla f'(\bar\rho_1)$ in the previous equation. After an integration
by parts, since $f'(\bar\rho_1) = \lambda_1$, in view of the boundary
conditions of \eqref{3-08b} and the definition \eqref{1-08xb} of
$\mf A_\lambda$, the right-hand side becomes
\begin{equation*}
\int_0^{\infty}\!dt \int_\Omega 
\big[ \, f'(u(t)) \,-\,  f'(\bar\rho_1) \, \big]
\; \nabla \cdot j(t) \;dx \;.
\end{equation*}
By \eqref{3-08b} this expression is equal to
\begin{equation*}
- \, \int_0^{\infty}\!dt \int_\Omega 
\big[ f'(u(t)) - f'(\bar\rho_1) \big] \, \partial_t u(t)  \; dx \;.
\end{equation*}

We have therefore shown that
\begin{equation}
\label{exw=qp}
W_\mathrm{ex} (E_1, \lambda_1, \bar\rho_0)  = \int_\Omega 
\big[ f(\bar\rho_0) - f(\bar\rho_1) -  
f'(\bar\rho_1) \big( \bar\rho_0 - \bar\rho_1\big) \big]
\; dx \;=\; V_{E_1, \lambda_1} (\bar\rho_0) \;,
\end{equation}
where the last identity follows from \eqref{x12}.

Note that the excess work $W_\mathrm{ex}$ is not the difference of a
thermodynamic potential between the states $\bar\rho_0$ and
$\bar\rho_1$.  We refer to \cite{bgjl1} and \cite[Ch.~7]{pippard} for
a connection of this result with \emph{availability} and the maximal
useful work that can be extracted from the system.

\section{Transformation along nonequilibrium states}
\label{sec5}

Nonequilibrium states are characterized by the presence of a non
vanishing current in the stationary density profile. Therefore, to
maintain such states one needs to dissipate a positive amount of
energy per unit of time.  If we consider a transformation between
nonequilibrium stationary states, the energy dissipated along such
transformation will necessarily include the contribution needed to
maintain such states. The arguments of the previous section have
therefore to be modified in order to take into account this amount of
energy.  This issue, first raised in \cite{op}, has been more recently
considered e.g.\ in \cite{bmw, maes2, kn, knst, bgjl1, bgjl2}.

The appropriate definition of thermodynamic functionals for
nonequilibrium systems is a central but difficult topic.  Our starting
point is the formula \eqref{04} for the energy exchanged in the time
interval $[0,T]$ between the system and the external reservoirs and
fields.

\subsection*{Towards a definition}

Fix $T>0$ and a nonequilibrium state $(E,\lambda)$ so that
$J_{E} (\bar\rho_{E,\lambda}) \not = 0$ and
$f'(\bar\rho_{E,\lambda}) \not = \lambda$. Let
$(E(t) , \lambda(t)) = (E,\lambda)$, $0\le t\le T$. By \eqref{04},
\begin{equation*}
\begin{aligned}
W_{[0,T]}(\lambda (\cdot), E (\cdot), \bar\rho_{E,\lambda}) 
\; & =\; T\, \int_\Omega 
J_{E} (\bar\rho_{E,\lambda}) \cdot \sigma(\bar\rho_{E,\lambda})^{-1}
J_{E} (\bar\rho_{E,\lambda})\; dx \\
& +\, T \;  \int_{\partial\Omega}
\mf A_\lambda (\, \bar\rho_{E,\lambda} \,,\,
f'(\bar\rho_{E,\lambda}) \,-\, \lambda \,) \; dS  \;.
\end{aligned}
\end{equation*}
Note that both terms on the right-hand side are strictly positive.
The second one is strictly positive in view of the last assertion of
\eqref{3-06} and because $f'(\bar\rho_{E,\lambda}) \,\not =\, \lambda$
at the boundary in nonequilibrium states.
\smallskip

To justify the definition of renormalized work proposed below, we turn
back to formula \eqref{04} for the work. Assume that the
transformation is performed along equilibrium states:
$J_{E(t)} (\bar\rho_{E(t),\lambda (t)}) =0$ for all $0\le t\le T$.

Since, in equilibrium, the current is equal to its symmetric part, in
\eqref{04}, $j(t) = J_{E(t)} (u(t)) = J^s_{E(t),\lambda (t)}
(u(t))$. On the other hand, in equilibrium, by equation \eqref{x10},
$(\delta V_{E,\lambda}/\delta \rho) (\rho) = f'(\rho) -
\lambda$. Therefore, for transformations along equilibrium states, we
may rewrite the work as
\begin{equation}
\label{5-01}
\begin{aligned}
& F (u(T)) \;-\;  F(\rho) 
\; +\; \int_{0}^{T}\! dt \int_\Omega
J^s_{E(t), \lambda (t)} (u(t)) 
\cdot \sigma(u(t))^{-1} J^s_{E(t), \lambda (t)} (u(t))  \; dx \\
& +\; \int_0^T \! dt  \int_{\partial\Omega}
\mf A_{\lambda(t)}  \Big(\, u(t) \,,\, 
\frac{\delta V_{E(t),\lambda (t)}}{\delta \rho} (u(t)) \,\Big)
\; dS \;.
\end{aligned}
\end{equation}

Fix an equilibrium state $(E,\lambda)$ and consider a constant
transformation $(E(t),\lambda(t)) = (E,\lambda)$, $0\le t\le T$,
starting from the density profile $\bar\rho_{E,\lambda}$.  Since
$\bar\rho_{E,\lambda}$ minimizes the quasi-potential $V_{E,\lambda}$,
$(\delta V_{E,\lambda}/\delta \rho) (\bar\rho_{E,\lambda}) =
0$. Hence, by the definition \eqref{1-08xb} of
$\mf A_{\lambda }$,
$\mf A_{\lambda } (\, \bar\rho_{E,\lambda} \,,\, (\delta V_{E
,\lambda}/\delta \rho)\, (\bar\rho_{E,\lambda}) \,) \,=\, \mf
H^{(3)}_{\lambda } (\, \bar\rho_{E,\lambda} \,,\, 0 \,) \,=\, 0$.  On
the other hand, by \eqref{3-25},
$J^s_{E, \lambda} (\bar\rho_{E,\lambda}) =0$. This shows that both
integrals in the previous displayed equation vanish.

This property extends to non-equilibrium states $(E,\lambda)$.
Indeed, since the quasi-potential $V_{E,\lambda}$ is minimal at the
stationary profile, $(\delta V_{E,\lambda}/\delta \rho)$
$(\bar\rho_{E,\lambda})=0$. Hence, by \eqref{3-25},
$J^s_{E,\lambda} (\bar\rho_{E,\lambda})=0$, and by the reasons
presented in the previous paragraph,
$\mf A_{\lambda } (\, \bar\rho_{E,\lambda} \,,\, 0 \,) \,=\, 0$.

\subsection*{Renormalized work}

The previous arguments support the following definition or
renormalized work.  Fix $T>0$, a density profile $\rho$, and
space-time dependent chemical potentials $\lambda(t)=\lambda(t,x)$ and
external field $E(t)=E(t,x)$, $0\le t\le T$, $x\in\Omega$.  Let
$u(t)=u(t, x)$, $j(t)=J_{E(t)} (u(t,x))$, $t \ge 0$, $x\in\Omega$, be
the solution of the hydrodynamic equation \eqref{3-08b} with initial
condition $\rho$.  Define the renormalized work
$W^\textrm{ren}_{[0,T]} = W^\textrm{ren}_{[0,T]}(E(\cdot),\lambda
(\cdot),\rho)$ performed by the reservoirs and the external field in
the time interval $[0,T]$ as
\begin{equation}
\label{5-02}
\begin{aligned}
W^\textrm{ren}_{[0,T]} \; & :=\; F (u(T)) \;-\;  F(\rho) \\
\;& +\; \int_{0}^{T}\! dt \int_\Omega
J^s_{E(t), \lambda (t)} (u(t)) 
\cdot \sigma(u(t))^{-1} J^s_{E(t), \lambda (t)} (u(t))  \; dx \\
& +\; \int_0^T \! dt  \int_{\partial\Omega}
\mf A_{\lambda(t)}  \Big(\, u(t) \,,\, 
\frac{\delta V_{E(t),\lambda (t)}}{\delta \rho} (u(t)) \,\Big)
\; dS \;.
\end{aligned}
\end{equation}

\begin{remark}
\label{rm6}
In \cite{bgjl1}, the renormalized work is defined by subtracting some
quantities from the work. Then, taking advantage of the orthogonality
between the symmetric and the anti-symmetric currents, a formula for
the renormalized work similar to \eqref{5-02} is derived.

Here, instead of subtracting a quantity we rather replaced in the
formula of the work the equilibrium quasi-potential by the
nonequilibrium one and the current by the symmetric current.

As observed in the previous subsection, the above definition of
renormalized work coincides with the one of work when the states
$(E(t), \lambda(t))$, $0\le t\le T$, are all equilibrium states.
\end{remark}

Assume that $\lambda (t), E(t)$ converge to $\lambda_1, E_1$ as
$t\to+\infty$ fast enough.  Let
$\bar\rho_1 = \bar \rho_{E_1,\lambda_1}$ be the stationary profile
associated to the pair $(E_1,\lambda_1)$.  Since $u(T)$ converges to
$\bar\rho_1$, the symmetric part of the current,
$J^s_{E(T), \lambda (T)}(u(T))$, relaxes as $T\to \infty$ to
$J^s_{E_1, \lambda_1}(\bar\rho_1) = 0$ fast enough.  Similarly,
$[\, \delta V_{E(T),\lambda (T)} / \delta \rho] \, (u(T))$ converges
to $[\, \delta V_{E_1,\lambda_1} / \delta \rho\,] \, (\bar\rho_1) =0$.
Hence, since, by \eqref{x2b}, $\mf A_{\lambda}$ is quadratic
in the second variable,
$\mf A_{\lambda(T)} (\, u(T) \,,\, [\, \delta V_{E(T),\lambda
(T)} / \delta \rho] \, (u(T)) \,)$ relaxes quickly to
$\mf A_{\lambda_1} (\, u(T)$,
$[\, \delta V_{E_1,\lambda_1} / \delta \rho\,] \, (\bar\rho_1) \,) =
\mf A_{\lambda_1} (\, u(T) \,,\, 0 \,) =0$.

Therefore, the two integrals in the previous formula are convergent as
$T\to \infty$ and
\begin{align}
W^\textrm{ren}_{[0,T]} \; & =\; F (\bar\rho_1) \;-\;  F(\rho)
\;+\; \int_{0}^{\infty}\! dt \, \int_{\partial\Omega}
\mf A_{\lambda(t)}  \Big(\, u(t) \,,\, 
\frac{\delta V_{E(t),\lambda (t)}}{\delta \rho} (u(t)) \,\Big)
\; dS
\nonumber\\
\;& +\; \int_{0}^{\infty}\! dt \int_\Omega
J^s_{E(t), \lambda (t)} (u(t)) 
\cdot \sigma(u(t))^{-1} J^s_{E(t), \lambda (t)} (u(t))  \; dx \;.
\label{5-04}
\end{align}
Since, by \eqref{x2b}, $\mf A_\lambda (a,p) \ge 0$,
\begin{equation}
\label{5-03}
W^\textrm{ren} (\, E(\cdot) \,,\, \lambda(\cdot) \,,\, \rho\, )
\;\ge\;  F (\bar\rho_1) \;-\;  F(\rho) \;.
\end{equation}
The previous equation states that the Clausius inequality holds for
the renormalized work, see \cite{kn, bgjl1}.

\subsection*{Quasi-static transformations}

As for transformations of equilibrium states, we show that, given two
nonequilibrium states, there exists a sequence of transformations from
the first to the second for which the last two terms on the right-hand
of \eqref{5-04} can be made arbitrarily small.

Fix $(E_0, \lambda_0)$ and assume that the initial profile $\rho$ is
the stationary profile associated to this pair:
$\rho = \bar\rho_{E_0, \lambda_0}$.  Fix $T>0$ and choose smooth
functions $\lambda(t),E(t)$, $0\le t\le T$, such that
$(E(0), \lambda(0))=(E_0, \lambda_0)$,
$(E(T), \lambda(T))=(E_1, \lambda_1)$.  For $\delta >0$, let
$(E_\delta(t), \lambda_\delta (t)) = (E(\delta t), \lambda (\delta
t))$, and $u_\delta(t)$ be the solution of \eqref{3-08b} with initial
condition $\bar\rho_0=\bar\rho_{E_0, \lambda_0}$, external field
$E_\delta(t)$, and chemical potential $\lambda_\delta (t)$.  Set
$j_\delta(t) = J_{E_\delta(t)}(u_\delta(t))$.  The last term on the
right-hand side of \eqref{5-04} is given by
\begin{equation*}
\int_{0}^{\infty}\!dt \int_\Omega
J^s_{E_\delta (t), \lambda_\delta (t)} (u_\delta(t))  \cdot 
\sigma(u_\delta(t))^{-1} J^s_{E_\delta (t),
\lambda_\delta (t)} (u_\delta(t)) \; dx \;.
\end{equation*}
For each fixed $t$, let
$\bar\rho_\delta(t) = \bar\rho_{E_\delta(t), \lambda_\delta(t)}$ be
the stationary profile associated to the driving
$E_\delta(t), \lambda_\delta (t)$ with frozen $t$.  Since
$J^s_{E_\delta (t), \lambda_\delta (t)} (\bar\rho_\delta(t)) =0$, we
can rewrite the previous integral as
\begin{equation*}
\int_{0}^{\infty}\!dt \int_\Omega
\widehat{J}_\delta (t) \cdot 
\sigma(u_\delta(t))^{-1} \widehat{J}_\delta (t) \; dx\;,
\end{equation*}
where
$\widehat{J}_\delta (t) = J^s_{E_\delta (t), \lambda_\delta (t)}
(u_\delta(t)) \,-\, J^s_{E_\delta (t), \lambda_\delta (t)}
(\bar\rho_\delta(t))$.

The difference between the solution of the hydrodynamic equation
$u_\delta(t)$ and the stationary profile $\bar\rho_\delta(t)$ is of
order $\delta$ uniformly in time, and so is the difference
$J^s_{E_\delta (t), \lambda_\delta (t)} (u_\delta(t)) \,-\,
J^s_{E_\delta (t), \lambda_\delta (t)} (\bar\rho_\delta(t))$.  As the
integration over time essentially extends over an interval of length
$\delta^{-1}$, the previous expression vanishes for $\delta\to 0$.

A similar argument can be carried out to the first integral in
\eqref{5-04} because $u_\delta(t)$ is close to $\bar\rho_\delta(t)$,
$[\, \delta V_{E_\delta (t), \lambda_\delta (t)} / \delta \rho\,] \,
(\bar\rho_\delta(t)) =0$,
$\mf A_{\lambda_\delta(t)} (\, u_\delta(t) \,,\, 0 \,) =0$, and
$\mf A_{\lambda}$ is quadratic in the second variable.  This
implies that equality in \eqref{5-03} is achieved in the limit
$\delta\to 0$.  In this argument we did not use any special property
of the path $(E(t), \lambda(t))$ besides its smoothness in time, the
trajectory $(E(t), \lambda (t))$ from $(E_0, \lambda_0)$ to
$(E_1, \lambda_1)$ can be otherwise arbitrary.

Quasi static transformations thus minimize asymptotically the
renormalized work and in the limit $\delta\to 0$ we obtain the
nonequilibrium version of the thermodynamic relation \eqref{12}, that is
\begin{equation}
\label{5-06}
W^\textrm{ren}  \;=\;  \Delta F \;,
\end{equation}
where $\Delta F$ represents the variation of the equilibrium free
energy functional, $\Delta F = F(\bar\rho_1) - F(\bar\rho_0)$.

It is remarkable that the Clausius inequality and the optimality of
quasi-static transformations, basic laws of equilibrium
thermodynamics, admit exactly the same formulation for nonequilibrium
states with the definition proposed in \eqref{5-02}. By Remark
\ref{rm6}, \eqref{5-03}, \eqref{5-06} contain the equilibrium
situations as a particular case.

\subsection*{Relaxation path: excess work and  quasi potential}

Consider at time $t=0$ a stationary nonequilibrium profile
$\bar\rho_0$ corresponding to some driving $(E_0, \lambda_0)$.  This
system is put in contact with new reservoirs at chemical potential
$\lambda_1$ and a new external field $E_1$.  For $t> 0$ the system
evolves according to the hydrodynamic equation \eqref{3-08b} with
initial condition $\bar\rho_0$, time independent boundary condition
$\lambda_1$ and external field $E_1$.  In particular, as $t\to \infty$
the system relaxes to $\bar\rho_1$.

Along such a path, in view of the orthogonality relation \eqref{3-15},
writing $J^s$ as $J-J^a$, the excess work is given by
\begin{equation*}
\begin{aligned}
W_\mathrm{ex} (E_1, \lambda_1,\bar\rho_0) 
\; & : = \;
\int_{0}^{\infty}\! dt \, \int_{\partial\Omega}
\mf A_{\lambda_1}  \Big(\, u(t) \,,\, 
\frac{\delta V_{E_1,\lambda_1}}{\delta \rho} (u(t)) \,\Big)
\; dS \\
\;& +\; \int_{0}^{\infty}\! dt \int_\Omega
J_{E_1} (u(t)) 
\cdot \sigma(u(t))^{-1} J^s_{E_1, \lambda_1} (u(t))  \; dx \\
& - \; \int_{0}^{\infty}\! dt  \int_{\partial \Omega}
\ms M^{\rm bd}_{\lambda_1, u(t)}
\big(\,  \frac{\delta V_{E_1,\lambda_1}}{\delta \rho} (u(t))\, \big)
\; \kappa \; d{\rm S} \;.
\end{aligned}
\end{equation*}

By definition \eqref{3-25} of the symmetric part of the current, by an
integration by parts, and by \eqref{3-08b}, \eqref{x8}, the second
term on the right-hand side is equal to
\begin{equation*}
\begin{aligned}
& \int_0^\infty \!dt \int_{\Omega} \nabla \cdot J_{E_1} (u(t)) 
\,  \frac{\delta V_{E_1, \lambda_1} }{\delta \rho} (u(t))  \; dx \\
& +\; \int_{0}^{\infty }\! dt \, \int_{\partial\Omega} 
\frac{\delta V_{E_1, \lambda_1}}{\delta \rho} (u(t))  \,
(\ms M^{\rm bd}_{\lambda_1, u(t)})' (0)\,
\kappa\; dS \;.
\end{aligned}
\end{equation*}
By \eqref{3-08b},
\begin{equation*}
\int_0^\infty \!dt \int_{\Omega} \nabla \cdot J_{E_1} (u(t)) 
\,  \frac{\delta V_{E_1, \lambda_1} }{\delta \rho} (u(t))  \; dx
\;=\; -\, \int_0^\infty \!dt \int_{\Omega} \partial_t u(t) 
\,  \frac{\delta V_{E_1, \lambda_1}}{\delta \rho} (u(t)) \; dx\;.
\end{equation*}
The right-hand side is equal to
$V_{E_1, \lambda_1} (\bar\rho_0) - V_{E_1, \lambda_1} (\bar\rho_1) =
V_{E_1, \lambda_1} (\bar\rho_0)$.

In conclusion, as $\ms M^{\rm bd}_{\lambda_1, u(t)} (0) =0$, in view
of the definition \eqref{1-08xb} of $\mf A_{\lambda_1}$,
\begin{align}
\nonumber
W_\mathrm{ex} (E_1, \lambda_1,\bar\rho_0) 
& \;=\; V_{E_1, \lambda_1} (\bar\rho_0) \;+\;
\int_{0}^{\infty}\! dt \, \int_{\partial\Omega}
\mf A_{\lambda_1}  \Big(\, u(t) \,,\, 
\frac{\delta V_{E_1,\lambda_1}}{\delta \rho} (u(t)) \,\Big)
\; dS \\
& - \; \int_{0}^{\infty}\! dt  \int_{\partial \Omega}
\mf A_{\lambda_1} \big(\, u(t) \, ,\,
\frac{\delta V_{E_1,\lambda_1}}{\delta \rho} (u(t))\, \big) \;
d{\rm S} \;,
\label{5-07}
\end{align}
so that
\begin{equation*}
W_\mathrm{ex} (E_1, \lambda_1,\bar\rho_0) 
\;=\; V_{E_1, \lambda_1} (\bar\rho_0)\;.
\end{equation*}
This identity extends to nonequilibrium states the relation \eqref{exw=qp}
between the excess work and the quasi potential.

\begin{remark}
As observed in \cite{bgjl1}, the previous identity provides a
characterization of the quasi-potential which does not involve large
deviations.
\end{remark}

\appendix

\section{Zero range dynamics}
\label{sec2}

Recall the notation introduced in Section \ref{sec1}.  In this
section, $\color{blue} \ms E = \bb N \cup \{0\}$ so that
$\color{blue} c(\ms E) = [0,\infty)$. The dynamics can be informally
described as follows.  At each site, independently from the others,
particles wait exponential times, whose parameter depends only on the
number of particles at that site, and then jumps to a nearest
neighboring site according to the transition probability of some
random walk on $\Omega_N$.  Superimposed to this bulk dynamics, to
model the effect of the reservoirs, we have creation and annihilation
of particles, according to some birth and death process, at the
boundary of $\Omega_N$.

Fix a time-dependent external field
$E\colon \bb R_+ \times \Omega \to \bb R^d$ and chemical potential
$\lambda \colon \bb R_+ \times \partial \Omega \to \bb R_+$.  The
generator $L_{t,N}$ of the zero range process is given by
\begin{equation*}
L_t \;=\;  L^{\rm bulk}_{t,N} \;+\;  L^{\rm bd}_{t,N} \;,
\end{equation*}
where $L^{\rm bulk}_{t,N}$ describes the bulk dynamics and
$ L^{\rm bd}_{t,N}$ the boundary dynamics at time $t$.  The generator
of the bulk dynamics is given by
\begin{equation*}
(L^{\rm bulk}_{t,N} f)(\eta) \;=\;
N^2 \, \sum_{x\in\Omega_N} \sum_{\substack{y\in\Omega_N\\ |y-x|=\varepsilon_N}} 
g(\eta_x) \, e^{  (1/2) E(t,x) \cdot (y-x)} \,
\big[ \, f(\sigma^{x,y}\eta) -f (\eta)\, \big] \;.
\end{equation*}
In this formula, $\color{blue} \varepsilon_N = 1/N$ and
$\sigma^{x,y}\eta$ is the configuration obtained from $\eta$ by moving
a particle from $x$ to $y$:
\begin{equation}
\label{9-01}
(\sigma^{x,y} \eta)_z  \;=\;
\left\{
\begin{array}{ccl}
\eta_z &\hbox{if}& z\neq x,y \\
\eta_z -1 &\hbox{if}& z=x \\
\eta_z +1  &\hbox{if}& z=y\;. 
\end{array}
\right.
\end{equation}

The generator of the boundary dynamics is given by
\begin{equation*}
\begin{aligned}
(L^{\rm bd}_{t,N} f)(\eta) \; & =\;
N \, \sum_{x\in\Omega_N} \sum_{\substack{y\not\in\Omega_N\\ |y-x|=\varepsilon_N}} 
g(\eta_x) \, e^{(1/2)  E(t,x) \cdot (y-x)} \,
\big[ \, f(\sigma^{x,-}\eta) -f (\eta)\, \big] \\
& + \;
N \, \sum_{x\in\Omega_N} \sum_{\substack{y\not\in\Omega_N\\ |y-x|=\varepsilon_N}} 
e^{\lambda (t,y)} \, e^{(1/2)  E(t,y) \cdot (x-y)} \,
\big[ \, f(\sigma^{x,+}\eta) -f (\eta)\, \big] \;.
\end{aligned}
\end{equation*}
In this formula, $\sigma^{x,+}\eta$, $\sigma^{x,-}\eta$ are the
configurations obtained from $\eta$ by removing, adding a particle at
$x$, respectively:
\begin{equation}
\label{9-02}
(\sigma^{x,\pm}\eta)_z
\left\{
\begin{array}{ccl}
\eta_z &\hbox{if}& z\neq x \\
\eta_z \pm 1 &\hbox{if}& z=x \;.
\end{array}
\right.
\end{equation}

Note that the bulk dynamics has been speeded-up by $N^2$, while the
boundary dynamics by $N$. Denote by $\color{blue} \eta^N(t)$ the
continuous-time Markov chain on $\Omega_N$ induced by the generator
$L_{t,N}$ and by $\color{blue} \bb P^{\lambda, E}_\eta$,
$\eta\in \Omega_N$, the distribution of the process $\eta^N(\cdot)$
when its initial state is $\eta$.

\subsection*{Stationary states}

Consider the case in which the driving $(\lambda, E)$ does not depend
on time. As the Markov chain is irreducible, there exists a
{\color{blue} unique invariant measure, denoted by
$\mu^{\lambda, E}_N$}.  It is remarkable that such invariant measure
can be constructed explicitly and it is product, see \cite{DF} for the
one dimensional case.

Denote by $\lambda_c \in \overline{\bb R}$ the radius of convergence
of the series
\begin{equation}
\label{Z=}
Z(\lambda) \;=\;  1 \;+\; \sum_{k\ge 1}
\frac {e^{\lambda k}}{g(1)
\cdots g(k)}
\end{equation}
For $\lambda < \lambda_c$, let $m_\lambda$ be the probability measure
on $\bb N$ given by
\begin{equation}
\label{2-04}
m_\lambda (k)  
\;=\; \frac {1}{Z(\lambda)} \; 
{\frac {e^{\lambda k}}{g(1)\cdots g(k)}} \;,
\quad k \in \bb N \cup \{0\} \;.
\end{equation}

Let $\phi_N = \phi^{\lambda, E}_N : \Omega_N \to \bb R_+$ be the
unique solution of the elliptic equation
\begin{equation*}
\begin{gathered}
\begin{aligned}
& N^2 \sum_{\substack{y\in\Omega_N\\ |y-x|=\varepsilon_N}}
\Big\{\, \phi_N (y) \,e^{  (1/2)  E(y) \cdot (x-y)} 
\,-\, \phi_N(x) \, e^{  (1/2)  E(x) \cdot (y-x)} \, 
\Big\} \\
+\; & N\,
\sum_{\substack{y\not\in\Omega_N\\ |y-x|=\varepsilon_N}}
\Big\{\, \phi_N (y) \,e^{  (1/2)  E(y) \cdot (x-y)} 
\,-\, \phi_N(x) \, e^{  (1/2)  E(x) \cdot (y-x)} \, 
\Big\} \;=\; 0\;, \quad x\in \Omega_N\;, \\
\end{aligned}
\\
\phi_N(z) \,=\, e^{\lambda(z)} \;, \quad
z\, \not\in\, \Omega_N \;.
\end{gathered}
\end{equation*}

The stationary state $\mu^{\lambda, E}_N$ is the product measure
on $\Omega_N$ whose marginals are given by
\begin{equation*}
\mu^{\lambda, E}_N \{\eta : \eta_x = k\}
\;=\; m_{\lambda_N(x)}(k) \;, \quad x\in \Omega_N\;,\;
k\ge 0\;, \;\;\text{where}\;\; \lambda_N(x) \;=\;
\log \phi_N (x)\;.
\end{equation*}

In the homogeneous equilibrium state, $E=0$ and $\lambda$ constant,
the solution of the elliptic equation is given by
$\phi_N=\exp\{\lambda\}$ so that the invariant measure is Gibbs with
Hamiltonian
\begin{equation*}
H_N (\eta) = \sum_{x\in\Omega_N} \sum_{k=1}^{\eta_x}  \log g(k)\;.
\end{equation*}

\subsection*{Bulk Hamiltonian}

Recall the definition of the function $R$ its inverse $\Xi = R^{-1}$,
and the equilibrium free energy, introduced in \eqref{3-2} and below.
By \cite{B8}, the diffusivity and the mobility are respectively given
by
\begin{equation*}
D(\rho) \;=\; \Phi'(\rho) \;, \quad
\sigma (\rho) \;=\; \Phi(\rho)  \;, \;\; \text{where}\;\;
\Phi(\rho)  \;=\;  e^{\Xi (\rho)} \;.
\end{equation*}
An elementary computation yields that the Einstein relation
\eqref{0-1} is fulfilled and that
$\color{blue} f'(\rho) = \log \Phi(\rho)$. By \eqref{3-00}, the bulk
Hamiltonian is given by
\begin{equation}
\label{3-00b}
\ms H^{\rm bulk}_{E} (\rho, F) \;=\; 
-\, \int_{\Omega} \Phi'(\rho) \, \nabla \rho \, \cdot \, \nabla F\; dx
\;+\; \int_{\Omega}  \Phi(\rho) \,
\big\{\, E \,+\, \nabla F\, \big\} \, \cdot \, \nabla F
\; dx  \;,
\end{equation}

\subsection*{Boundary Hamiltonian}

In view of the definition of the generator $L^{\rm bd}_{t,N}$, in the
context of the zero-range process, the boundary generator
$\mc L_\lambda$, $\lambda<\lambda_c$, introduced in \eqref{3-1}, is
given by
\begin{equation*}
(\mc L_{\lambda} f)(\mf x) \;=\; g(\mf x)\, [\, f(\mf x-1) \,-\, f(\mf
x)\,] \;+\; e^\lambda \, [\, f(\mf x+1) \,-\, f(\mf x)\,]   \;.
\end{equation*}

The boundary Hamiltonian,
$\ms M^{\rm bd}_{\lambda} \colon (0,\infty)\times \bb R \to \bb R$,
introduced in \eqref{3-02}, is given by
\begin{equation}
\label{3-04b}
\ms M^{\rm bd}_{\lambda, a}(p) \;=\;
e^\lambda \, [e^p-1] \;+\; e^{\Xi(a)} \, [e^{-p}-1]  \;.
\end{equation}

\subsection*{Quasi-potential}

Assume that the external field vanishes and fix a chemical potential
$\lambda$. Recall the definition of the variable $d(\cdot)$ introduced
below \eqref{3-13}. In the context of zero-range processes,
$\color{blue} d(\rho) = \Phi(\rho) = \sigma(\rho)$.

We claim that
\begin{equation}
\label{x20}
\frac{\delta V_{0,\lambda}}{\delta \rho} (\rho) \;=\;
f'(\rho) \;-\; f'(\bar\rho_{0,\lambda}) \;,
\end{equation}
where $\bar\rho_{0,\lambda}$ the solution of \eqref{3-09}.  In
particular, for zero-range processes, the quasi-potential has an
explicit formula.

To prove \eqref{x20}, we first claim that identity \eqref{x18} holds
for zero-range processes. Since $d(\rho) = \sigma(\rho) = \Phi(\rho)$,
and $f'(\rho) = \log \Phi(\rho)$, the left-hand side of \eqref{x18}
can be written as
\begin{equation*}
\frac{\Phi(\rho)}{\Phi(\rho) - \Phi(F)}\,
\Big\{\, e^\lambda \, \Big[\, \frac{\Phi(\rho)}{\Phi(F)} \,-\,
1\,\Big]
\;+\; \Phi(\rho) \, \Big[\, \frac{\Phi(F)}{\Phi(\rho)} \,-\,
1\,\Big] \,\Big\} \;=\; e^\lambda \;-\; \Phi(F)\;.
\end{equation*}
The right-hand side is equal to $(\ms M^{\rm bd}_{\lambda, F})'(0)$,
proving  \eqref{x18}.

Since $d(\rho) = \sigma(\rho)$, by Remark \ref{rm-x} and
\eqref{3-131}, in the context of zero-range processes, equation
\eqref{3-13} becomes
\begin{equation*}
\left\{
\begin{aligned}
& \Delta f'(F) \;+\; 
\, \Vert \nabla f'(F)\, \Vert^2 \;=\; 0 \;, \\
& J(F) \cdot \bs n \;=\; -\,
\kappa\, (\ms M^{\rm bd}_{\lambda, F})' (0) \;.
\end{aligned}
\right.
\end{equation*}
A simple algebra based on Einstein relation and the relations between
mobility, diffusivity and $\Phi$ permits to rewrite the previous
equation as
\begin{equation}
\label{3-13A}
\left\{
\begin{aligned}
& \nabla \cdot J(F)  \;=\; 0 \;, \\
& J(F) \cdot \bs n \;=\; -\,
\kappa\, (\ms M^{\rm bd}_{\lambda, F})' (0) \;.
\end{aligned}
\right.
\end{equation}
This equation corresponds to the stationary equation \eqref{3-09}.
Hence, by \eqref{3-14}, for zero-range dynamics, the quasi-potential
is given by \eqref{x20}.

\section{Exclusion processes}
\label{sap2}

Recall the notation introduced in Section \ref{sec1}. In the context
of exclusion processes, $\ms E = \{0,1\}$ so that $c(\ms E) = [0,1]$.
The dynamics can be informally described as follows. Particles are
distributed on $\Omega_N$ in such a way that, at each site, there is
at most one particle.  Each particle, independently from the others,
wait a mean-one exponential random time, and then jumps to a nearest
neighboring site according to the transition probability of some
random walk on $\Omega_N$. If the chosen site is occupied by another
particle, the jumps is suppressed. Superimposed to this bulk dynamics,
to model the effect of the reservoir, at the boundary of $\Omega_N$,
particles are created and annihilated according to some birth and
death process.

Fix a time-dependent external field
$E\colon \bb R_+ \times \Omega \to \bb R^d$ and chemical potential
$\lambda \colon \bb R_+ \times \partial \Omega \to \bb R_+$.  The
generator $L_{t,N}$ of the exclusion process is given by
\begin{equation}
\label{3-23}
L_{t,N} \;=\;  L^{\rm bulk}_{t,N} \;+\;  L^{\rm bd}_{t,N} \;,
\end{equation}
where $L^{\rm bulk}_{t,N}$ describes the bulk dynamics and
$ L^{\rm bd}_{t,N}$ the boundary dynamics at time $t$.  The generator
of the bulk dynamics is given by
\begin{equation*}
(L^{\rm bulk}_{t,N} f)(\eta) \;=\;
N^2 \, \sum_{x\in\Omega_N} \sum_{\substack{y\in\Omega_N\\ |y-x|=\varepsilon_N}} 
\eta_x \, [\, 1-\eta_y\,] \, e^{  (1/2) E(t,x) \cdot (y-x)} \,
\big[ \, f(\sigma^{x,y}\eta) -f (\eta)\, \big] \;,
\end{equation*}
where $\varepsilon_N = 1/N$ and $\sigma^{x,y}\eta$ have been introduced in
\eqref{9-01}. 

The generator of the boundary dynamics is given by
\begin{equation*}
\begin{aligned}
& (L^{\rm bd}_{t,N} f)(\eta) \; =\;
N \, \sum_{x\in\Omega_N} \sum_{\substack{y\not\in\Omega_N\\ |y-x|=\varepsilon_N}} 
\eta_x \, \frac{1}{1+ e^{\lambda (t,y)}}\, e^{(1/2)  E(t,x) \cdot (y-x)} \,
\big[ \, f(\sigma^{x,-}\eta) -f (\eta)\, \big] \\
& \quad + \;
N \, \sum_{x\in\Omega_N} \sum_{\substack{y\not\in\Omega_N\\ |y-x|=\varepsilon_N}} 
[\, 1 \,-\, \eta_x \,]\,
\frac{e^{\lambda (t,y)}}{1+ e^{\lambda (t,y)}}  \, e^{(1/2)  E(t,y) \cdot (x-y)} \,
\big[ \, f(\sigma^{x,+}\eta) -f (\eta)\, \big] \;,
\end{aligned}
\end{equation*}
where the configuration $\sigma^{x,\pm}\eta$ has been introduced in
\eqref{9-02}. 

Note that the bulk dynamics has been speeded-up by $N^2$, while the
boundary dynamics by $N$. Denote by $\color{blue} \eta^N(t)$ the
continuous-time Markov chain on $\Omega_N$ induced by the generator
$L_{t,N}$ and by $\color{blue} \bb P^{\lambda, E}_\eta$,
$\eta\in \Omega_N$, the distribution of the process $\eta^N(\cdot)$
when its initial state is $\eta$.

\subsection*{Stationary states}

Consider the case in which the driving $(\lambda, E)$ does not depend
on time. As the Markov chain is irreducible, there exists a
{\color{blue} unique invariant measure, denoted by
$\mu^{\lambda, E}_N$}. In contrast with the zero-range process, beyond
the equilibrium case where the current vanishes, the stationary state
$\mu^{\lambda, E}_N$ is not a product measure and exhibits long range
correlations \cite{Spo}. In the special case $E=0$, $\lambda$, the
measure $\mu^{\lambda, E}_N$ is the product measure with Bernoulli
marginals of density $e^\lambda/(1+e^\lambda)$.

\subsection*{Bulk Hamiltonian}

By \cite{B8}, the diffusivity and the mobility are respectively given
by
\begin{equation*}
D(\rho) \;=\; 1 \;, \quad
\sigma (\rho) \;=\; \rho\, (1-\rho) \;.
\end{equation*}
An elementary computation yields that the Einstein relation
\eqref{0-1} is fulfilled. By \eqref{3-00}, the bulk Hamiltonian is
given by
\begin{equation}
\label{3-00c}
\ms H^{\rm bulk}_{E} (\rho, F) \;=\; 
-\, \int_{\Omega}  \nabla \rho \, \cdot \, \nabla F\; dx
\;+\; \int_{\Omega}  \rho\, (1-\rho) \,
\big\{\, E \,+\, \nabla F\, \big\} \, \cdot \, \nabla F
\; dx  \;,
\end{equation}

\subsection*{Boundary Hamiltonian}

In view of the definition of the generator $L^{\rm bd}_{t,N}$, in the
context of the exclusion process, the boundary generator
$\mc L_\lambda$, introduced in \eqref{3-1}, is given by
\begin{equation*}
(\mc L_{\lambda} f)(0) \;=\; \frac{e^\lambda}{1+e^\lambda}\,
[\, f(1) \,-\, f(0)\,] \;, \quad
(\mc L_{\lambda} f)(1) \;=\;
\frac{1}{1+e^\lambda}  \, [\, f(0) \,-\, f(1)\,]   \;.
\end{equation*}

The boundary Hamiltonian,
$\ms M^{\rm bd}_{\lambda} \colon [0,1]\times \bb R \to \bb R$,
introduced in \eqref{3-02}, is given by
\begin{equation}
\label{3-03c}
\ms M^{\rm bd}_{\lambda, \rho}(p) \;=\; [1-\rho]\, R(\lambda) \,
[e^p-1] \;+\; \rho\, [\,1-R(\lambda)\,] \, [e^{-p}-1]  \;,
\end{equation}
where $R(\lambda)$ is the mean of the measure $m_\lambda$ and has been
introduced in \eqref{3-2}.

\subsection*{Quasi-potential}

Assume that $d=1$, $\Omega = (0,1)$, and recall the definition of the
function $d(\cdot)$ introduced below \eqref{3-13}. For exclusion
processes, $\color{blue} d(\rho) = \rho$,
$f'(\rho) = \log [\rho /(1-\rho)]$. By \eqref{3-03c}, equation
\eqref{x18} is satisfied, and equation \eqref{3-131} takes the form
\begin{equation}
\label{3-13c}
\left\{
\begin{aligned}
& \Delta F  \;=\; (\, \rho \,-\, F\, )\, 
\, \frac{(\, \nabla F \,)^2}{F\,(1-F)}  \;, \\
& F' \cdot \bs n \,=\, \kappa \, [\, \varrho (\lambda) - F \,]
\;\; \text{at}\;\;  x=0\;, \;\; x=1\;,
\end{aligned}
\right.
\end{equation}
where $\varrho (\lambda) = e^\lambda /(1+e^\lambda)$.  It has been
shown in \cite{BEL21} that equation \eqref{3-13c} has a unique
solution. Hence, by \eqref{3-14}, in dimension $1$ with no external
field, 
\begin{equation*}
\frac{\delta V_{0,\lambda}}{\delta \rho} (\rho) \;=\;
f'(\rho) \;-\; f'(F) \;.
\end{equation*}

\section{KMP model}
\label{sap3}

Recall the notation introduced in Section \ref{sec1}. This time
$\ms E = c(\ms E) = \bb R_+$, and $\eta_x$, $x\in\Omega_N$, represents
the energy at site $x$ for the configuration $\eta$.  The bulk
dynamics can be informally described as follows. At each bond $(x,y)$
in $\Omega_N$, at exponential times, the energy of the two vertices is
added and then redistributed according to a uniform measure.

Fix a time-dependent chemical potential
$\lambda \colon \bb R_+ \times \bb R^d \to \bb R_-$. \emph{Note that
$\lambda$ takes negative values}. We adopted this convention, which
might be slightly confusing, to uniformize the notation of all three
models. Moreover, there is no external field. The generator $L_{t,N}$
of the KMP process is given by
\begin{equation*}
L_t \;=\;  L^{\rm bulk}_{t,N} \;+\;  L^{\rm bd}_{t,N} \;,
\end{equation*}
where $L^{\rm bulk}_{t,N}$ describes the bulk dynamics and
$ L^{\rm bd}_{t,N}$ the boundary dynamics at time $t$.  The generator
of the bulk dynamics is given by
\begin{equation*}
(L^{\rm bulk}_{t,N} f)(\eta) \;=\;
N^2 \, \sum_{(x,y)\in\Omega_N}
\int_0^1 
\big[ \, f(\sigma^{x,y}_r\eta) - f (\eta)\, \big] \; dr \;.
\end{equation*}
In this formula, the sum is performed over all unordered edges of
$\Omega_N$ and $\sigma^{x,y}_r\eta$ is the configuration obtained
from $\eta$ by replacing $\eta_x$, $\eta_y$ by $r(\eta_x+\eta_y)$,
$(1-r)(\eta_x+\eta_y)$, respectively:
\begin{equation*}
(\sigma^{x,y}_r \eta)_z  \;=\;
\left\{
\begin{array}{ccl}
\eta_z &\hbox{if}& z\neq x,y \\
r(\eta_x+\eta_y) &\hbox{if}& z=x \\
(1-r)(\eta_x+\eta_y)  &\hbox{if}& z=y\;. 
\end{array}
\right.
\end{equation*}

The generator at the boundary is given by
\begin{equation*}
(L^{\rm bd}_{t,N} f)(\eta) \;=\;
N \, \sum_{x\in\Omega_N}
\sum_{\substack{y\not\in\Omega_N\\ |y-x|=\varepsilon_N}} 
\int_0^\infty -\, \lambda(t,y) \, e^{\lambda (t,y) r}
\big[ \, f(\sigma^{x}_r\eta) - f (\eta)\, \big] \; dr \;,
\end{equation*}
where $\sigma^{x}_r\eta$ is the configuration obtained
from $\eta$ by replacing $\eta_x$ by $r$:
\begin{equation*}
(\sigma^{x}_r \eta)_z  \;=\;
\left\{
\begin{array}{ccl}
\eta_z &\hbox{if}& z\neq x \\
r &\hbox{if}& z=x \;. 
\end{array}
\right.
\end{equation*}

\subsection*{The stationary states}

Denote by $m_\lambda$, $\lambda <0$, the distribution of an
exponential random variable with mean
$\color{blue} \tau(\lambda) = - \lambda^{-1}$.  Denote by
$\mu^{\lambda}_N$ the product measure on $\Sigma_N$ whose marginals
are given by
\begin{equation*}
\int_{\Sigma_N} F(\eta_x) \; \mu^{\lambda}_N (d\eta)
\;=\; \int_{\bb R_+} F(\mf x) \;
m_{\lambda}(d\mf x) \;, \quad x\in \Omega_N\;,\;\; F\in C_b(\bb R)\;. 
\end{equation*}
An elementary computation shows that $\mu^{\lambda}_N$ is a stationary
state (actually, reversible) for the KMP dynamics when $\lambda(t,x)$
is constant and equal to $\lambda$. Ergodicity yields that it is the
unique one.

\subsection*{Bulk Hamiltonian}

By \cite{B8}, the diffusivity and the mobility of the KMP models are
given by
\begin{equation}
\label{10-01}
D(\rho) \;=\; 1\;, \quad \sigma(\rho) \;=\; \rho^2\;.
\end{equation}
Therefore, by \eqref{3-00}, the bulk Hamiltonian is given by
\begin{equation}
\label{3-00d}
\ms H^{\rm bulk}(\rho, F) \;=\; 
-\, \int_{\Omega} \nabla \rho \, \cdot \, \nabla F\; dx
\;+\; \int_{\Omega}  \rho^2 \,
\nabla F\, \cdot \, \nabla F \; dx  \;.
\end{equation}

\subsection*{Boundary Hamiltonian}

For KMP dynamics, the boundary generator $\mc L_\lambda$,
$\lambda < 0$, introduced in \eqref{3-1}, is given by
\begin{equation*}
(\mc L_{\lambda} f)(\mf x) \;=\;
(-\, \lambda)\, \int_{\bb R_+} [\, f(\mf y) \,-\, f(\mf
x)\,] \,  e^{\lambda \mf y} \, d\mf y  \;.
\end{equation*}

The boundary Hamiltonian,
$\ms M^{\rm bd}_{\lambda, \rho} \colon (0,\tau^{-1}) \to \bb R$,
introduced in \eqref{3-02}, is given by
\begin{equation}
\label{3-05b}
\ms M^{\rm bd}_{\lambda, \rho}(p) \;=\;
\frac{\tau}{\rho+\tau} \,\Big(\frac{1}{1-\tau \, p}\,-\, 1\Big) \;+\;
\frac{\rho}{\rho+\tau} \,\Big(\frac{1}{1+\rho \, p}\,-\, 1\Big) \;.
\end{equation}

\subsection*{Quasi-potential}

Assume that $d=1$, $\Omega = (0,1)$, and recall the definition of the
function $d(\cdot)$ introduced below \eqref{3-13}. For KMP dynamics,
$\color{blue} d(\rho) = \rho$, $f'(\rho) = - (1/\rho)$, and equation
\eqref{3-13} takes the form
\begin{equation}
\label{3-13d}
\left\{
\begin{aligned}
& \Delta F  \;+\; (\, \rho \,-\, F\, )\, 
\, \frac{(\, \nabla F\,)^2}{F^2}  \;=\; 0 \;, \\
& \nabla F  \cdot \bs n \;=\;
\kappa\, F^2\, \frac{\tau - F}{\rho F - \tau \rho + \tau F}
\;.
\end{aligned}
\right.
\end{equation}

In terms of the variables $G=f'(F)$, $\gamma = f'(\rho)$, $\lambda =
f'(\tau)$, the equation reads
\begin{equation}
\label{3-13db}
\left\{
\begin{aligned}
& \Delta G  \;-\; \Big(\, \frac{1}{\gamma} + \frac{1}{G}\, \Big)\, 
\, (\nabla G)^2   \;=\; 0 \;, \\
& \nabla G  \cdot \bs n \;=\;
\kappa\, \gamma\, \frac{G-\lambda}{G - \lambda - \gamma}
\;.
\end{aligned}
\right.
\end{equation}

Uniqueness of solutions of equation \eqref{3-13d} has still to be
proven.  It has been done in \cite{BGL05} if the boundary conditions
are replaced by Dirichlet ones.  If uniqueness holds, by \eqref{3-14},
the quasi-potential of the KMP model is given by
\begin{equation*}
\frac{\delta V_{0,\lambda}}{\delta \rho} (\rho) \;=\;
f'(\rho) \;-\; f'(F) \;=\; \frac{1}{F} \;-\; \frac{1}{\rho} \;\cdot
\end{equation*}

\section{Exclusion process with non-reversible boundary conditions}
\label{sap4}

Inspired by the model introduced in \cite{DPTV11, DPTV12}, and the
works \cite{ELX18, EGN20}, in this section, we present a model which
do not satisfy \eqref{x13}, \eqref{3-06}. This is a consequence from
the fact that the stationary state induced by the boundary dynamics
does not coincide with the one induced by the bulk dynamics. To
concentrate on the source of the differences, we assume that there is
no external field and we set $d=1$.

Let $\Omega =(0,1)$ so that
$\Omega_N = \{\epsilon_N , \dots, 1-\epsilon_N\}$, where, recall,
$\epsilon_N = 1/N$.  The state space and the bulk dynamics are the
ones introduced in Section \ref{sap2} with $\Omega = (0,1)$ and
$E=0$. To define the boundary dynamics we introduce a set of jump
rates. Fix $\ell\ge 1$, and let
$c^R_j \colon \{0,1\}^{\{-\ell, \dots, -1\}} \to \bb R_+$,
$c^L_j \colon \{0,1\}^{\{1, \dots, \ell\}} \to \bb R_+$,
$1\le j\le \ell$ be nonnegative functions.

The generator of the boundary dynamics is given by
\begin{equation*}
(L^{\rm bd}_{t,N} f)(\eta) \; =\;
N \, \sum_{j=1}^{\ell} 
c^R_{j}(\tau_{N} \eta)\,
\big[ \, f(\sigma^{N-j}\eta) -f (\eta)\, \big]
\;+\;
N \, \sum_{j=1}^{\ell} 
c^L_{j}(\eta)\,
\big[ \, f(\sigma^{j}\eta) -f (\eta)\, \big]\;.
\end{equation*}
In this formula, $\color{blue} \tau_N \eta$ is the configuration
$\eta$ translated by $N$ so that $(\tau_N \eta)_j = \eta_{N+j}$,
$j\in \bb Z$.  Moreover, the configuration $\sigma^{k}\eta$ stands for
\begin{equation*}
(\sigma^{k}\eta)_j
\left\{
\begin{array}{ccl}
\eta_j &\hbox{if}& j\neq k \\
1-\eta_k  &\hbox{if}& j=k \;.
\end{array}
\right.
\end{equation*}
The generator $L_N$ of the dynamics is given by \eqref{3-23}, and does not
depend on time.

The model introduced above embraces the exclusion process introduced
in Section \ref{sap2} with no external field (to incorporate the
chemical potential, it is enough to let the jump rates $c^{R,L}_j$ to
depend on $\lambda$). It also encompasses the current reservoir model
considered by De Masi et al. \cite{DPTV11, DPTV12} and the exclusion
models with nonreversible boundary dynamics examined in \cite{ELX18,
EGN20}.

\subsection*{Bulk $\times$ Boundary dynamics}

In contrast with the previous models, here particles are created and
annihilated at more than one site in the bulk of $\Omega_N$, and
according to different rates which depend on the environment. For this
reason, the boundary dynamics can not be represented by a one-site
dynamics as in \eqref{3-1}.  The state space is here
$\{0,1\}^{\{1, \dots, \ell\}}$ instead of $\{0,1\}$ as in the
exclusion dynamics of Section \ref{sap2}.

Consider a neighborhood $\{N-k, \dots, N-1\}$ of the right boundary.
The stationary state of the bulk dynamics restricted to this set (we
forbid exchange of particles between $N-k-1$ an $N-k$) is
the uniform measure over all configurations with a fixed number of
particles. For $k$ large, by the equivalence of ensembles, locally
this measure is close to a Bernoulli product measure with some fixed
density.

Unless in very special cases, the stationary state on
$\{N-\ell, \dots, N-1\}$ induced by the generator $L^{\rm bd}_{t,N}$
introduced above is not a Bernoulli product measure. When this does
not happen, there is a conflict between the bulk dynamics, which
drives the system towards a Bernoulli product measure, and the bulk
dynamics, which propels the system to another stationary state. As the
bulk dynamics is accelerated by $N^2$, while the boundary dynamics is
speeded-up by $N$, the bulk dynamics wins and the state of the system
at the boundary is close the a Bernoulli product measure. In
particular, local equilibrium occurs and the entropy method to derive
the hydrodynamic behavior can be applied \cite{EGN20}.

One can also, up to technical obstacles, prove a large deviations
principle and derive a formula for the Hamiltonian. A rigorous proof
of this statement is not yet available.

\subsection*{Stationary states}

Assume that the Markov chain induced by the generator $L_N$ is
irreducible. This is the case, for example, if one of the jump rates
of each side of the set $\Omega_N$ is strictly positive.  In this
case, there exists a unique invariant measure, denoted by
$\mu_N$. Except in exceptional cases, the stationary state is not a
product measure and not known explicitly.

\subsection*{Bulk Hamiltonian}

The bulk Hamiltonian is the one presented in \eqref{3-00c}.

\subsection*{Boundary Hamiltonian}

Since the boundary of $\Omega$ consists of two points,
$\partial \Omega = \{0, 1\}$, the surface integral becomes a sum and
the boundary Hamiltonian reads
\begin{equation*}
\ms H^{\rm bd} \big(\, \rho\, ,\, F \, \big)  \;:=\;
\ms M^{\rm bd, 0}_{\rho}  (\, F \,) \; \kappa(0)
\;+\; \ms M^{\rm bd, 1}_{\rho}  (\, F \,) \; \kappa(1)\;, 
\end{equation*}
where
\begin{gather*}
\ms M^{\rm bd, 1}_{\rho}  (p) \;=\; \sum_{j=1}^{\ell}
E_{\nu_\rho} \Big[\,   c^R_{j}(\eta)\,\,
\big [\, e^{p(1-2\eta_{-j})} -1 \, \big] \,\Big] \;, \\
\ms M^{\rm bd, 0}_{\rho}  (p) \;=\; \sum_{j=1}^{\ell}
E_{\nu_\rho} \Big[\,   c^L_{j}(\eta)\,\,
\big [\, e^{p(1-2\eta_{j})} -1 \, \big] \,\Big] \;.
\end{gather*}
In this formula, $\nu_\rho$ represents the Bernoulli product measure
with density $\rho$.  These expressions can be written as in
\eqref{3-03c}. For $k=0$, $1$,
\begin{equation}
\label{3-03d}
\ms M^{\rm bd, k}_{\rho}(p) \;=\; [1-\rho]\, R^+_k(\rho) \,
[e^p-1] \;+\; \rho\, R^-_k (\rho) \, [e^{-p}-1]  \;,
\end{equation}
where
\begin{equation*}
R^+_1(\rho) \;=\; \sum_{j=1}^{\ell}
E_{\nu_\rho} \big[\,   c^R_{+,j}(\eta)\, \,\big]\;, \quad
R^-_1(\rho) \;=\; \sum_{j=1}^{\ell}
E_{\nu_\rho} \big[\,   c^R_{-,j}(\eta)\, \,\big]\;.
\end{equation*}
Here, $ c^R_{\pm,j}(\eta) =  c^R_{j}(\eta^{\pm,-j})$, and $\eta^{\pm,-j}$ is
the configuration which coincides with $\eta$ at all sites but $-j$,
and at $-j$ takes the value $[1 - (\pm 1)]/2$:
$(\eta^{-,-j})_{-j}=1$ and $(\eta^{+,-j})_{-j}=0$. A similar formula
holds for $R^\pm_0(\rho)$.

\subsection*{What does not hold for this model }

The proof of \eqref{x13} presented in Section \ref{sec1} requires the
boundary dynamics to be stationary with respect the measure induced by
the bulk dynamics, a property which does not hold here. In
consequence, the arguments presented in Section \ref{sec0} to show
that in equilibrium the boundary density satisfies the identity
$f'(\rho) = \lambda$ do not apply. Relation \eqref{x13} is also used
below \eqref{x2b} to rewrite the boundary term in \eqref{x7} as an
integral of the functional $\mf A_\lambda$. In particular, the proof
of Clausius inequality does not apply to this model.

This means that few assertions made in this article remain valid for
this model which deserves further investigations.

\subsection*{Quasi-potential}

An explicit formula for the quasi-potential, similar to \eqref{3-13c},
is an open problem for this model.

\begin{remark}
The equilibrium stationary states of the microscopic dynamics
presented in these appendices are all product measures. It should be
possible to extend this theory to one-dimensional Ising models under
the Kawasaki dynamics in mild contact with boundary reservoirs. The
technical difficulty lies in the fact that these models are
non-gradient \cite[Chapter 7]{kl}. 
\end{remark}

\smallskip\noindent{\bf Acknowledgments.} The authors wish to thank
the referees for their careful reading of a previous version of this
article and D. Gabrielli for many insights.

C. L. has been partially supported by FAPERJ CNE E-26/201.207/2014, by
CNPq Bolsa de Produtividade em Pesquisa PQ 303538/2014-7, by
ANR-15-CE40-0020-01 LSD of the French National Research Agency.


\begin{thebibliography}{99}


\bibitem{A89} V. I. Arnol'd: {\it Mathematical methods of classical
mechanics}, second edition, Springer-Verlag, New York,  1989.

\bibitem{bmw} M. Baiesi, C. Maes, B. Wynants: Fluctuations and
response of nonequilibrium states.  Phys. Rev. Lett. \textbf{103},
010602 (2009).

\bibitem{bmns} R. Baldasso, O. Menezes, A.  Neumann, R. R. Souza:
Exclusion process with slow boundary. J. Stat. Phys. {\bf 167},
1112--1142 (2017).

\bibitem{B2} L. Bertini, A. De Sole, D. Gabrielli, G. Jona-Lasinio,
and C. Landim, Macroscopic fluctuation theory for stationary non
equilibrium state, J. Statist. Phys. {\bf 107}, 635--675 (2002).

\bibitem{B8} L. Bertini, A. De Sole, D. Gabrielli, G. Jona-Lasinio, C.
Landim; Large deviation approach to non equilibrium processes in
stochastic lattice gases. Bol. Soc. Brasil. Mat. (N.S.)  {\bf 37}, 611
--643 (2006).

\bibitem{bdgjl11} L. Bertini, A. De Sole, D. Gabrielli,
G. Jona-Lasinio, C.  Landim; Towards a nonequilibrium thermodynamics:
a self-contained macroscopic description of driven diffusive systems.
J. Stat. Phys. {\bf 135}, 857--872 (2009).

\bibitem{bdgjl14} L. Bertini, A. De Sole, D. Gabrielli,
G. Jona-Lasinio, C.  Landim; Macroscopic fluctuation
theory. Rev. Modern Phys. {\bf 87}, 593--636 (2015).

\bibitem{bgjl1} L. Bertini, D. Gabrielli, G. Jona-Lasinio, and
C. Landim: Thermodynamic transformations of nonequilibrium states.
J. Statist. Phys. {\bf 149}, 773--802 (2012).

\bibitem{bgjl2} L. Bertini, D. Gabrielli, G. Jona-Lasinio, C.  Landim;
Clausius inequality and optimality of quasi static transformations for
nonequilibrium stationary states.  Phys. Rev. Lett. {\bf 110}, 020601
(2013)

\bibitem{BGL05} L. Bertini, D. Gabrielli, J. L. Lebowitz:
Large deviations for a stochastic model of heat flow.
J. Statist. Phys. {\bf 121}, 843--885, (2005).

\bibitem{maes2} E. Boksenbojm, C. Maes, K. Neto\u{c}n\'y,
J. Pe\u{s}ek: Heat capacity in nonequilibrium steady states.
Europhys. Lett. EPL \textbf{96}, 40001 (2011).

\bibitem{BEL21} A. Bouley, C. Erignoux, C. Landim: Steady state large
deviations for one-dimensional, symmetric exclusion processes in weak
contact with reservoirs. arXiv:2107.06606 (2021)

\bibitem{DF} De Masi A., Ferrari P.; A remark on the hydrodynamics of
the zero-range processes.  J. Stat. Phys. \textbf{36}, 81–87
(1984).

\bibitem{DPTV11} A. De Masi, E. Presutti, D. Tsagkarogiannis,
M. E. Vares: Current Reservoirs in the Simple Exclusion Process.  J
Stat Phys {\bf 144}, 1151--1170 (2011).

\bibitem{DPTV12} A. De Masi, E. Presutti, D. Tsagkarogiannis,
M. E. Vares: Non-equilibrium Stationary States in the Symmetric Simple
Exclusion with Births and Deaths. J Stat Phys {\bf 147}, 519--528
(2012). https://doi.org/10.1007/s10955-012-0481-2

\bibitem{DHS} B. Derrida, O. Hirschberg, T Sadhu: Large deviations in
the symmetric simple exclusion process with slow
boundaries. J. Stat. Phys. (2021) doi.org/10.1007/s10955-020-02680-3

\bibitem{DLS} B. Derrida, J. L. Lebowitz, E. R. Speer: Large
deviation of the density profile in the steady state of the open
symmetric simple exclusion process. J. Stat. Phys. {\bf 107}, 599--634
(2002).

\bibitem{EGN20} C. Erignoux, P. Gon\c{c}alves, G. Nahum:
Hydrodynamics for SSEP with non-reversible slow boundary dynamics:
Part I, the critical regime and beyond. J Stat Phys {\bf 181},
1433--1469 (2020). https://doi.org/10.1007/s10955-020-02633-w

\bibitem{ELX18} C. Erignoux, C. Landim, T. Xu: Stationary states of
boundary driven exclusion processes with nonreversible boundary
dynamics.  J. Stat. Phys. {\bf 171}, 599--631 (2018).

\bibitem{FGLN2021} T. Franco, P. Gon{\c c}alves, C. Landim, A.
Neumann: Dynamical large deviations for boundary driven symmetric
exclusion processes in weak contact with reservoirs. arXiv:2203.14417
(2022)

\bibitem{fv17} S. Friedli, Y. Velenik: {\it Statistical Mechanics of
Lattice Systems: A Concrete Mathematical Introduction} Cambridge
University Press, 2017. 

\bibitem{kl} C. Kipnis, C. Landim \emph{Scaling limits of interacting
particle systems.}  Springer-Verlag, Berlin, 1999.

\bibitem{kn} T. Komatsu, N. Nakagawa: Expression for the stationary
distribution in nonequilibrium steady states.
Phys. Rev. Lett. \textbf{100}, 030601 (2008).

\bibitem{knst} T. Komatsu, N. Nakagawa, S. Sasa, H. Tasaki: Entropy
and nonlinear nonequilibrium thermodynamic relation for heat
conducting steady states.  J. Stat. Phys. \textbf{142}, 127--153
(2011).

\bibitem{op} Y. Oono, M. Paniconi: Steady state thermodynamics.
Dynamic organization of fluctuations (Nishinomiya, 1997).
Progr. Theoret. Phys. Suppl. \textbf{130}, 29–44 (1998).

\bibitem{pippard} A. B.Pippard: \emph{Elements of classical
thermodynamics for advanced students of physics.}  Cambridge
University Press, New York 1957.

\bibitem{Spo} H. Spohn: Long range correlations for stochastic lattice
gases in a nonequilibrium steady state.  J. Phys. A {\bf 16},
4275--4291 (1983).

\bibitem{v} S. R. S. Varadhan: {\it Large deviations and
applications}. Society for Industrial and Applied Mathematics, 1984.

\end{thebibliography}
\end{document}